\documentclass{amsart}

\usepackage{amsthm}
\usepackage{amsmath}
\usepackage{amssymb}
\usepackage{amsfonts}
\usepackage{latexsym}
\usepackage{enumitem}
\usepackage{mathrsfs}
\usepackage[all]{xy} \SelectTips{eu}{}
\usepackage{hyperref}

\newtheorem{intthm}{Theorem}[]

\newcommand{\numberseries}{\bfseries}   %Fontseries used for numbering
\newlength{\thmtopspace}                %Space above theorem
\newlength{\thmbotspace}                %Space below theorem
\newlength{\thmheadspace}               %Space after theorem label
\newlength{\thmindent}                  %For indenting
\newtheoremstyle{bfupright head,slanted body}
                {\thmtopspace}{\thmbotspace}
                {\slshape}{\thmindent}{\bfseries}{.}{\thmheadspace}
                {{\numberseries \thmnumber{#2\;}}\thmnote{#3}}

\newtheoremstyle{bfupright head,upright body}
                {\thmtopspace}{\thmbotspace}
                {\upshape}{\thmindent}{\bfseries}{.}{\thmheadspace}
                {{\numberseries \thmnumber{#2\;}}\thmnote{#3}}

\newtheoremstyle{fixed bf head,slanted body}
                {\thmtopspace}{\thmbotspace}{\slshape}
                {\thmindent}{\bfseries}{.}{\thmheadspace}
                {{\numberseries \thmnumber{#2\;}}\thmname{#1}\thmnote{ (#3)}}

\newtheoremstyle{fixed bf head,upright body}
                {\thmtopspace}{\thmbotspace}{\upshape}
                {\thmindent}{\bfseries}{.}{\thmheadspace}
                {{\numberseries \thmnumber{#2\;}}\thmname{#1}\thmnote{ (#3)}}

\newtheoremstyle{numbered paragraph}
                {\thmtopspace}{\thmbotspace}{\upshape}
                {\thmindent}{\upshape}{}{\thmheadspace}
                {{\numberseries \thmnumber{#2.}}}

\theoremstyle{bfupright head,slanted body}
\newtheorem{res}{}[section]             \newtheorem*{res*}{}

\theoremstyle{bfupright head,upright body}
\newtheorem{bfhpg}[res]{}               \newtheorem*{bfhpg*}{}

\theoremstyle{fixed bf head,slanted body}
\newtheorem{thm}[res]{Theorem}          \newtheorem*{thm*}{Theorem}
\newtheorem{prp}[res]{Proposition}      \newtheorem*{prp*}{Proposition}
\newtheorem{cor}[res]{Corollary}        \newtheorem*{cor*}{Corollary}
\newtheorem{lem}[res]{Lemma}            \newtheorem*{lem*}{Lemma}

\theoremstyle{fixed bf head,upright body}
\newtheorem{dfn}[res]{Definition}       \newtheorem*{dfn*}{Definition}
\newtheorem{rmk}[res]{Remark}           \newtheorem*{rmk*}{Remark}
           \newtheorem*{exa*}{Example}
\theoremstyle{numbered paragraph}
\newtheorem{ipg}[res]{}

\newlength{\thmlistleft}        %leftmargin
\newlength{\thmlistright}       %rightmargin
\newlength{\thmlistpartopsep}   %partopsep
\newlength{\thmlisttopsep}      %topsep
\newlength{\thmlistparsep}      %parsep
\newlength{\thmlistitemsep}     %itemsep

\newcounter{eqc}
  {\end{list}}%

\newcounter{prt}
\newenvironment{prt}{\begin{list}{\upshape (\alph{prt})}%
    {\usecounter{prt}%
      \setlength{\leftmargin}{\thmlistleft}%
      \setlength{\labelwidth}{\thmlistleft}%
      \setlength{\rightmargin}{\thmlistright}%
      \setlength{\partopsep}{\thmlistpartopsep}%
      \setlength{\topsep}{\thmlisttopsep}%
      \setlength{\parsep}{\thmlistparsep}%
      \setlength{\itemsep}{\thmlistitemsep}}}%
  {\end{list}}%

\newcounter{rqm}
\newenvironment{rqm}{\begin{list}{\upshape (\arabic{rqm})}%
    {\usecounter{rqm}%
      \setlength{\leftmargin}{\thmlistleft}%
      \setlength{\labelwidth}{\thmlistleft}%
      \setlength{\rightmargin}{\thmlistright}%
      \setlength{\partopsep}{\thmlistpartopsep}%
      \setlength{\topsep}{\thmlisttopsep}%
      \setlength{\parsep}{\thmlistparsep}%
      \setlength{\itemsep}{\thmlistitemsep}}}%
  {\end{list}}%

\newenvironment{itemlist}{\nopagebreak \begin{list}{$\bullet$}%
    {\setlength{\leftmargin}{1.5em}%
      \setlength{\labelwidth}{\thmlistleft}%
      \setlength{\rightmargin}{\thmlistright}%
      \setlength{\partopsep}{\thmlistpartopsep}%
      \setlength{\topsep}{\thmlisttopsep}%
      \setlength{\parsep}{\thmlistparsep}%
      \setlength{\itemsep}{\thmlistitemsep}}}%
  {\end{list}}%

\newenvironment{prf*}[1][Proof]{%
  \begin{proof}[\bf #1]
    \setcounter{equation}{0}
    }
  {\end{proof}
}

\newcommand{\pgref}[1]{\ref{#1}}

\renewcommand{\eqref}[1]{(\pgref{eq:#1})}

\newcommand{\thmcite}[2][?]{\cite[thm.~#1]{#2}}
\newcommand{\prpcite}[2][?]{\cite[prop.~#1]{#2}}
\newcommand{\corcite}[2][?]{\cite[cor.~#1]{#2}}
\newcommand{\lemcite}[2][?]{\cite[lem.~#1]{#2}}

\newcommand{\dfncite}[2][?]{\cite[def.~#1]{#2}}

\numberwithin{equation}{res}

\def\urltilda{\kern -.15em\lower .7ex\hbox{\~{}}\kern .04em}

\newcommand{\GF}[1]{\mathsf{GF}(#1)}
\newcommand{\PGF}[1]{\mathsf{PGF}(#1)}
\newcommand{\GI}[1]{\mathsf{GI}(#1)}
\newcommand{\F}[1]{\mathsf{Flat}(#1)}
\newcommand{\W}{\mathsf{W}}
\newcommand{\A}{\mathsf{A}}
\newcommand{\T}{\mathsf{T}}
\newcommand{\Cot}[1]{\mathsf{Cot}(#1)}

\newcommand{\Mod}[1]{\mathsf{Mod}(#1)}
\newcommand{\Prj}[1]{\mathsf{Prj}(#1)}

\newcommand{\xra}[2][]{\xrightarrow[#1]{\:#2\:}}

\newcommand{\ZZ}{\mathbb{Z}}

\newcommand{\Rop}{R^{\sf op}}
\newcommand{\Qop}{Q^{\sf op}}
\newcommand{\aop}{a^{\sf op}}

\newcommand{\colim}{\mbox{\rm colim}}

\newcommand{\C}{\mbox{\rm C}}

\newcommand{\Ker}[1]{\nobreak{\operatorname{Ker}#1}}
\newcommand{\Coker}[1]{\nobreak{\operatorname{Coker}#1}}

\newcommand{\Hom}[3][]{\operatorname{Hom}_{#1}(#2,#3)}
\newcommand{\tp}[3][Q]{\nobreak{#2\otimes_{#1}#3}}
\newcommand{\Ext}[4][R]{\operatorname{Ext}_{#1}^{#2}(#3,#4)}
\newcommand{\Rep}[2]{{\sf Rep}(#1,#2)}
\newcommand{\is}{\cong}
\newcommand{\dis}{\:\is\:}
\def\bF{{\boldsymbol F}}
\def\bE{{\boldsymbol E}}

   \def\soft#1{\leavevmode\setbox0=\hbox{h}\dimen7=\ht0\advance
    \dimen7 by-1ex\relax\if t#1\relax\rlap{\raise.6\dimen7
    \hbox{\kern.3ex\char'47}}#1\relax\else\if T#1\relax
    \rlap{\raise.5\dimen7\hbox{\kern1.3ex\char'47}}#1\relax
    \else\if d#1\relax\rlap{\raise.5\dimen7\hbox{\kern.9ex
    \char'47}}#1\relax\else\if D#1\relax\rlap{\raise.5\dimen7
    \hbox{\kern1.4ex\char'47}}#1\relax\else\if l#1\relax
    \rlap{\raise.5\dimen7\hbox{\kern.4ex\char'47}}#1\relax
    \else\if L#1\relax\rlap{\raise.5\dimen7\hbox{\kern.7ex
    \char'47}}#1\relax\else\message{accent \string\soft
    \space #1 not defined!}#1\relax\fi\fi\fi\fi\fi\fi}

\begin{document}

\title[Gorenstein flat representations of left rooted quivers]%
{Gorenstein flat representations of left rooted quivers}

\author[Z.\ Di]{Zhenxing Di}

\address{Zhenxing Di: \ Northwest Normal University, Lanzhou 730070, China}

\email{dizhenxing@163.com}

\author[S.\ Estrada]{Sergio Estrada}

\address{Sergio Estrada: \ Universidad de Murcia, Murcia 30100, Spain}

\email{sestrada@um.es}

\urladdr{http://webs.um.es/sestrada}

\author[L.\ Liang]{Li Liang}

\address{Li Liang (corresponding author): \ Lanzhou Jiaotong University, Lanzhou 730070, China}

\email{lliangnju@gmail.com}

\urladdr{https://sites.google.com/site/lliangnju}

\author[S.\ Odaba\c{s}\i]{Sinem Odaba\c{s}\i}

\address{Sinem Odaba\c{s}\i: \ Universidad Austral de Chile, Valdivia, Chile}

\email{sinem.odabasi@uach.cl}

\thanks{Z. Di was partly supported by NSF of China grant 11971388; S. Estrada was partly supported by grant MTM2016-77445-P (AEI/FEDER,UE) and Fundaci\'on Seneca grant 19880/GERM/15; L. Liang was partly supported by NSF of China grant 11761045 and the Foundation of A Hundred Youth Talents Training Program of Lanzhou Jiaotong University; S. Odaba\c{s}\i\ was partly supported by the research grant CONICYT/FONDECYT/Iniciaci\'{o}n/11170394.}

\date{\today}

\keywords{Left rooted quiver, Gorenstein flat representation, model structure}

\subjclass[2010]{18G25; 16G20}

\begin{abstract}
We study Gorenstein flat objects in the category $\Rep{Q}{R}$ of representations of a left rooted quiver $Q$
with values in $\Mod{R}$, the category of all left $R$-modules,
where $R$ is an arbitrary associative ring.
We show that a representation $X$ in $\Rep{Q}{R}$ is Gorenstein flat if and only if for each vertex $i$
the canonical homomorphism $\varphi_i^X: \oplus_{a:j\to i}X(j)\to X(i)$ is injective,
and the left $R$-modules $X(i)$ and $\Coker\,\varphi_i^X$ are Gorenstein flat.
As an application of this result,
we show that there is a hereditary abelian model structure on $\Rep{Q}{R}$
whose cofibrant objects are precisely the Gorenstein flat representations,
fibrant objects are precisely the cotorsion representations, and
trivial objects are precisely the representations
with values in the right orthogonal category of all projectively coresolved Gorenstein flat left $R$-modules.
\end{abstract}

\maketitle

\thispagestyle{empty}
\section*{Introduction}%%%%%%%%%%%%%%%%%%%%%%%%%%%%%%%%%%
\noindent
The study of the module-valued representations of quivers is an important topic
in current research in representation theory of groups and algebras.
Our present work is motivated by a series of results on module-valued representations over left rooted quivers.
Before introducing the results that inspired us, let us introduce some necessary notation.
Throughout this section, $Q=(Q_0,Q_1)$ is a left rooted quiver,
where $Q_0$ and $Q_1$ are the sets of all vertices and arrows in $Q$, respectively.
Denote by $\Rep{Q}{R}$ the category of representations of $Q$ with values in $\Mod{R}$.
For every vertex $i\in Q_0$ and every representation $X$ in $\Rep{Q}{R}$,
there exists a canonical homomorphism $$\varphi_i^X: \oplus_{a\in Q_{1}^{\ast\to i}}X(s(a))\to X(i),$$
where $s(a)$ is as usual the source of the arrow $a\in Q_1$.

Enochs and Estrada gave in \thmcite[3.1]{EE05} the characterization of a projective representation
by showing that a representation $X$ in $\Rep{Q}{R}$ is projective if and only if $\varphi_i^X$ is a monomorphism,
and the left $R$-modules $X(i)$ and $\Coker\,\varphi_i^X$ are projective for each vertex $i\in Q_0$.
A similar characterization for flat representations (colimits of projective representations) was proved by Enochs, Oyonarte and Torrecillas in \thmcite[3.7]{EOT04}.
They are exactly the representations $X$ in $\Rep{Q}{R}$
for which $\varphi_i^X$ is a monomorphism, and the $R$-modules $X(i)$ and $\Coker\,\varphi_i^X$ are flat for all vertex $i\in Q_0$.

As two primary generalizations of modules having Gorenstein dimension 0 introduced by Auslander and Bridger \cite{AusBri} to modules that are not
finitely generated, Gorenstein projective modules and Gorenstein flat modules were introduced by Enochs, Jenda and Torrecillas in \cite{EEnOJn95b} and \cite{EJT-93}, and further treated by Holm in \cite{HHl04a}.
Eshraghi, Hafezi and Salarian gave in \thmcite[3.5.1]{EHS13}
the structure of a Gorenstein projective representation in $\Rep{Q}{R}$.
Surprisingly, the characterization of Gorenstein projective representations
carries over the form of the above projective ones without more assumptions.
They showed that a representation $X$ in $\Rep{Q}{R}$ is Gorenstein projective
if and only if $\varphi_i^X$ is a monomorphism,
and the left $R$-modules $X(i)$ and $\Coker\,\varphi_i^X$ are Gorenstein projective for all vertex $i\in Q_0$.

In the present paper, we study Gorenstein flat representations in $\Rep{Q}{R}$.
The first step in describing Gorenstein flat representations is to
identify the notion of tensor products in $\Rep{Q}{R}$.
Let $X$ be a representation in $\Rep{Q}{R}$.
Inspired by the work of Salarian and Vahed \cite{SV2016},
we construct a tensor product functor $\tp{-}{X}$ (see \ref{tproduct}), which has enough nice expected attributes. For example, a representation $X$ in $\Rep{Q}{R}$ is flat if and only if the functor $\tp{-}{X}$ is exact; see Theorem \ref{flat}. With the notion of tensor products in hand,
we can give the definition of Gorenstein flat representations in $\Rep{Q}{R}$ in the routine way; see Definition \ref{gfqr}.
Our first main result asserts that a Gorenstein flat representation in $\Rep{Q}{R}$ can be described similarly as follows.

\begin{intthm}\label{AA}
Let $Q$ be a left rooted quiver and $X\in\Rep{Q}{R}$.
Then $X$ is Gorenstein flat if and only if $X$ is in $\Phi(\GF{R})$, that is, for each vertex $i\in Q_0$ the homomorphism
$\varphi_i^X: \oplus_{a\in Q_{1}^{\ast\to i}}X(s(a))\to X(i)$ is injective,
and the left $R$-modules $X(i)$ and $\C_i(X)$ are Gorenstein flat.
\end{intthm}

\begin{center}
  $\ast \ \ \ast \ \ \ast$
\end{center}
\noindent
Hovey studied extensively in \cite{Ho02} model structures on abelian categories.
The most celebrated result in \cite{Ho02}, which is now known as {\sf Hovey's correspondence},
says that an abelian model structure on an abelian category $\sf{A}$
is equivalent to a triple of subcategories ($\sf{C,W,F}$) in $\sf{A}$
such that $\sf{W}$ is thick and
$(\sf{C},\sf{W}\cap\sf{F})$
and $(\sf{C}\cap\sf{W},\sf{F})$ form two complete cotorsion pairs
(here, $\sf{W}$ (resp., $\sf{C}$ and $\sf{F}$) is the subcategory of $\A$
consisting of all trivial (resp., cofibrant and fibrant) objects
associated to the corresponding abelian model structure).
Hovey's correspondence makes it clear that an abelian model structure on $\sf{A}$
can be succinctly represented by the triple $\sf{(C,W,F)}$.
Therefore, one often refers to such a triple as an abelian model structure in literatures,
and call it a {\sf Hovey triple}.

Holm and J\o rgensen \cite{HJ19} proved that
if an abelian category $\sf{A}$ satisfies some mild conditions,
then every complete hereditary cotorsion pair in $\sf{A}$ induces two hereditary cotorsion pairs
in $\Rep{Q}{\sf{A}}$,
which gives a ``quiver representation version" of \corcite[3.8]{Gi04}
by Gillespie on cotorsion pairs in the categories of chain complexes.
Odaba\c{s}{\i} further showed in \cite{Od19} that the two induced cotorsion pairs are complete as well.
These wonderful results help us to obtain the following result,
which asserts that a hereditary Hovey triple $\sf{(C,W,F)}$ on $\sf{A}$
can induce a hereditary Hovey triple on $\Rep{Q}{\sf{A}}$
with trivial objects are precisely the representations with values in $\sf{W}$
(here by ``hereditary Hovey triple"
we mean that the two associated cotorsion pairs in Hovey's correspondence are hereditary); see Section 1 for unexplained notation.

\begin{intthm}\label{BB}
Let $Q$ be a left rooted quiver and
$\sf{A}$ an abelian category with enough injectives
which satisfies the axiom $\sf{AB4}$ (e.g. Grothendieck category).
Then any hereditary Hovey triple $\sf{(C,W,F)}$ on $\sf{A}$
induces an hereditary Hovey triple
$$(\Phi({\sf C}),\Rep{Q}{\W},\Rep{Q}{\sf{F}})$$
on $\Rep{Q}{\sf{A}}$.
\end{intthm}

Gillespie \cite{Gil17} constructed, over any coherent ring $\Lambda$, a new hereditary abelian model structure,
the Gorenstein flat model structure, on $\Mod{\Lambda}$
by showing that there exists a thick subcategory $\mathsf{W}(\Lambda)$
such that $(\GF{\Lambda}, \mathsf{W}(\Lambda), \Cot{\Lambda})$ forms a hereditary Hovey triple on $\Mod{\Lambda}$
(here $\GF{\Lambda}$ and $\Cot{\Lambda}$ denote the subcategory of $\Mod{\Lambda}$
consisting of all Gorenstein flat and cotorsion $\Lambda$-modules, respectively).
Recently, \v{S}aroch and \v{S}t'ov\'{\i}\v{c}ek \cite{SS20} extended the above Gorenstein flat model structure to an arbitrary ring $R$; they showed that
$(\GF{R},\PGF{R}^\perp,\Cot{R})$ forms a hereditary Hovey triple on $\Mod{R}$,
where $\PGF{R}^\perp$ is the right orthogonal category of all projectively coresolved Gorenstein flat $R$-modules.
As an application of Theorems \ref{AA} and \ref{BB},
we obtain a Gorenstein flat model structure on $\Rep{Q}{R}$ in which we give explicit descriptions of the subcategories of trivial, cofibrant and fibrant objects.

\begin{intthm}
Let $Q$ be a left rooted quiver.
Then there is a hereditary Hovey triple
$$(\GF{Q}, \PGF{Q}^\perp, \Cot{Q})$$
on $\Rep{Q}{R}$
in which $\GF{Q}=\Phi(\GF{R})$, $\PGF{Q}^\perp=\Rep{Q}{\PGF{R}^\perp}$ and $\Cot{Q}=\Rep{Q}{\Cot{R}}$.
\end{intthm}

The paper is organized as follows.
Section 1 contains some necessary notation and terminology for use throughout this paper.
In Section 2, we give the definition of tensor product functors in $\Rep{Q}{R}$,
and show the relation with flat representations.
In Section 3, we introduce the notion of Gorenstein flat representations in $\Rep{Q}{R}$,
and give the proof of Theorem \ref{AA}.
Section 4 is devoted to giving the proof of Theorem \ref{BB}.
Finally, we give an Appendix for reproving two key Lemmas \ref{cogenerator} and \ref{trivial objects}
over a more simple quiver with 4 vertices
to comprehend the ideas of the proofs.

\section{Preliminaries}\label{pre}%%%%%%%%%%%%%%%%%%%%%%%%%%%%%%%%%%%
\noindent
Throughout the paper,
all rings are assumed to be associative with identity.
Let $R$ be a ring;
we adopt the convention that an $R$-module is a left $R$-module,
and we refer to right $R$-modules as modules over the opposite ring $\Rop$.
We let $\Mod{R}$ denote the category of all $R$-modules,
and denote by $\GF{R}$ (resp., $\F{R}$, $\Prj{R}$ and $\Cot{R}$)
the subcategory of $\Mod{R}$ consisting of all Gorenstein flat (resp., flat, projective and cotorsion) $R$-modules.

\begin{bfhpg}[\bf Gorenstein flat/injective modules]
An $R$-module $M$ is called \emph{Gorenstein flat} (see e.g.  Enochs and Jenda \cite{rha}) if there is an exact sequence
$$\cdots \to F^{-1} \to F^0 \to F^{1} \to \cdots$$ of flat $R$-modules with $M\is\Ker{(F^0 \to F^1)}$, such that it remains exact after applying the functor $\tp[R]{E'}{-}$ for every injective $\Rop$-module $E'$.
The subcategory of all Gorenstein flat $R$-modules is denoted $\GF{R}$.

An $R$-module $M$ is called \emph{Gorenstein injective} (see e.g. \cite{rha}) if there is an exact sequence
$$\cdots \to I^{-1} \to I^0 \to I^{1} \to \cdots$$
of injective $R$-modules with $M\is\Ker{(I^0 \to I^1)}$, such that it remains exact after applying the functor $\Hom[R]{E}{-}$ for every injective $R$-module $E$.
The subcategory of all Gorenstein injective $R$-modules is denoted $\GI{R}$.
\end{bfhpg}

\begin{bfhpg}[\bf Cotorsion pairs]\label{cotpair}
Let $\sf A$ be an abelian category with enough projectives and injectives.
A pair $(\sf X,\sf Y)$ of subcategories of $\sf A$ is called a \emph{cotorsion pair} if $\sf{X}^{\bot}=\sf{Y}$ and $^{\bot}\sf{Y}=\sf{X}$. Here
$${\sf X}^{\bot}=\{N\in \sf A\ |\ \Ext[\sf A]{1}{X}{N}=0 {\rm \ for\ all}\ X\in\sf X\},\ \text{and}$$
$$^{\bot}{\sf Y}=\{M\in \sf A\ |\ \Ext[\sf A]{1}{M}{Y}=0 {\rm \ for\ all}\ Y\in\sf Y\}.$$
Recall from \cite{rha} that a cotorsion pair $(\sf X,\sf Y)$ is \emph{complete} if for each object $M$ in $\sf A$ there is an exact squence $0 \to M \to Y \to X \to 0$ with $X\in\sf X$ and $Y\in\sf Y$ (or equivalently, there is an exact squence $0 \to Y' \to X' \to M \to 0$ with $X'\in\sf X$ and $Y'\in\sf Y$). A cotorsion pair $(\mathsf{X},\mathsf{Y})$ is called \emph{hereditary} if $\Ext[\sf A]{i}{X}{Y}=0$ for each $X\in\mathsf{X}$ and $Y\in\mathsf{Y}$, and all $i\geq1$. A cotorsion pair $(\mathsf{X},\mathsf{Y})$ is called \emph{perfect} if all objects in $\sf A$ have $\sf X$-covers and $\sf Y$-envelopes in the sense of \cite{rha}.
\end{bfhpg}

\begin{bfhpg}[\bf Quivers and representations]\label{quiver notation1}
A quiver $Q$ is actually a directed graph with vertex set $Q_0$ and arrow set $Q_1$. The symbol $\Qop$ denotes the opposite quiver of $Q$, which is obtained by reversing arrows of $Q$. For an arrow $a\in Q_1$, we always write $s(a)$ for its source and $t(a)$ for its target. For a vertex $i\in Q_0$ we let $Q_1^{i\to\ast}$ denote the set $Q_1^{i\to\ast}=\{a\in Q_1~|~s(a)=i\}$, and let $Q_1^{\ast\to i}$ denote the set $Q_1^{\ast\to i}=\{a\in Q_1~|~t(a)=i\}$.

Let $Q$ be a quiver and $\A$ an abelian category.
By a \emph{representation} $X$ of $Q$ we mean a functor from $Q$ to $\A$,
which is determined by associating an object in $\A$ to each vertex $i\in Q_0$ and
a morphism $X(a): X(i)\to X(j)$ to each arrow $a: i\to j$ in $Q_1$.

Let $X$ and $Y$ be representations of $Q$.
A \emph{morphism} $f$ from $X$ to $Y$ is a natural transformation,
that is, a family of morphisms $\{f(i):X(i)\to Y(i)\}_{i\in Q_0}$
such that $Y(a)\circ f(i)=f(j)\circ X(a)$ for each arrow $a: i\to j$.

The representations of $Q$ constitute a Grothendieck category with enough projectives,
which is denoted $\Rep{Q}{\A}$.
In particular, if $\A=\Mod{R}$ then we use the notation $\Rep{Q}{R}$ instead of $\Rep{Q}{\Mod{R}}$.
\end{bfhpg}

\begin{ipg}\label{quiver notation}
Let $Q$ be a quiver.
For each $X\in\Rep{Q}{\A}$ and each vertex $i\in Q_0$,
by the universal property of coproducts,
there is a unique morphism $$\varphi_{i}^{X}:\oplus_{a\in Q_{1}^{\ast\to i}}X(s(a))\to X(i).$$
Let $\C_i(X)$ denote the cokernel of $\varphi_{i}^{X}$;
it yields a functor from $\Rep{Q}{\A}$ to $\A$.
Dually, there is a unique morphism
$$\psi_{i}^{X}:X(i)\to \Pi_{a\in Q_{1}^{i\to \ast}}X(t(a)).$$
The symbol $\mathrm{K}_i(X)$ denotes the kernel of $\psi_{i}^{X}$,
which also yields a functor from $\Rep{Q}{\A}$ to $\A$.

For a subcategory $\sf X$ of $\A$, we set
\begin{itemlist}
\item $\Rep{Q}{{\sf X}}=\{X\in\Rep{Q}{\A}~|~X(i)\in\sf X\}$,
\item $\Phi({\sf X})=\left\{ X\in\Rep{Q}{\A} \:
      \left|
        \begin{array}{c}
          \varphi_i^X\ \text {is a monomorphism and}\\
          \C_i(X),X(i)\in{\sf X}\ \text{for all}\ i\in Q_0
        \end{array}
      \right.
    \right\}$ and
\item $\Psi({\sf X})=\left\{ X\in\Rep{Q}{\A} \:
      \left|
        \begin{array}{c}
          \psi_i^X\ \text {is an epimorphism and}\\
          \mathrm{K}_i(X),X(i)\in{\sf X}\ \text{for all}\ i\in Q_0
        \end{array}
      \right.
    \right\}$
\end{itemlist}
\end{ipg}

\begin{bfhpg}[\bf Left rooted quivers]\label{rooted}
Let $Q$ be a quiver. It follows from \cite{EOT04}
that there is a transfinite sequence $\{V_{\alpha}\}_{\alpha \mathrm{ordinal}}$ of subsets of $Q_0$ as follows:

For the first ordinal $\alpha=0$ set $V_0=\emptyset$, for a successor ordinal $\alpha+1$ set
$$V_{\alpha+1}=\{i\in Q_0~|~i\ \text{is not the target of any arrow}\ a\in Q_1\ \text{with}\ s(a)\notin\cup_{\beta\leq\alpha}V_{\beta}\},$$
and for a limit ordinal $\alpha$ set $V_{\alpha}=\cup_{\beta<\alpha}V_{\beta}$.

It is clear that $V_1=\{i\in Q_0~|~\text{there is no arrow}\ a\in Q_1\ \text{with}\ t(a)=i\}$. By \lemcite[2.7]{HJ19} and \corcite[2.8]{HJ19}, there is a chain
$$V_1\subseteq V_2\subseteq\cdots\subseteq Q_0,$$ and
if $a: i\to j$ is an arrow in $Q_1$ with $j\in V_{\alpha+1}$ for some ordinal $\alpha$, then $i$ must be in $V_{\alpha}$.

Following \dfncite[3.5]{EOT04} a quiver $Q$ is called \emph{left rooted}
if there exists an ordinal $\lambda$ such that $V_{\lambda}=Q_0$.
By \prpcite[3.6]{EOT04}, a quiver $Q$ is left rooted if and only if
it has no infinite sequence of arrows of the form $\cdots \to \bullet\to \bullet \to \bullet$ (not necessarily different).
So the left rooted quivers constitute quite a large class of quivers.
\end{bfhpg}

The next result is concluded by \thmcite[7.4, 7.6 and 7.9]{HJ19} and \thmcite[4.6]{Od19}.
Recall that an abelian category $\sf{A}$ satisfies the axiom $\sf{AB4}$
provided that $\sf{A}$ is cocomplete such that any coproduct of monomorphisms in $\sf{A}$ is a monomorphism.

\begin{lem}\label{induced cp}
Let $Q$ be a left rooted quiver and $\sf{A}$ an abelian category with enough injectives
which satisfies the axiom $\sf{AB4}$ (e.g. Grothendieck category).
If $(\sf X, \sf Y)$ is a complete hereditary cotorsion pair in $\A$,
then $(\Phi(\sf X),\Rep{Q}{\sf Y})$ is a complete hereditary cotorsion pair in $\Rep{Q}{\A}$.
\end{lem}

\begin{bfhpg}[\bf Flat and cotorsion representations]
Recall from \cite{EOT04} that a representation $X\in\Rep{Q}{R}$ is \emph{flat}
if it is a colimit of projective representations.
The subcategory of all flat representations in $\Rep{Q}{R}$ is denoted $\F{Q}$.
A representation $X$ in $\Rep{Q}{R}$ is called \emph{cotorsion}
if it is in $\F{Q}^{\perp}$.
The subcategory of all cotorsion representations in $\Rep{Q}{R}$ is denoted  $\Cot{Q}$.
\end{bfhpg}

The next result can be found in \thmcite[3.7]{EOT04}.

\begin{lem}\label{flat}
Let $Q$ be a left rooted quiver and $X\in\Rep{Q}{R}$. Then $X$ is flat if and only if $X$ is in $\Phi(\F{R})$.
\end{lem}

The equality in the statement of the following result was first proved by Oyonarte \thmcite[6]{Oy01}.

\begin{lem}\label{cotorsion}
Let $Q$ be a left rooted quiver. Then $(\F{Q},\Cot{Q})$ is a complete hereditary cotorsion pair with $\Cot{Q}=\Rep{Q}{\Cot{R}}$.
\end{lem}
\begin{prf*}
It follows from \prpcite[2]{BEE-01} that $(\F{R},\Cot{R})$ is a complete hereditary cotorsion pair in $\Mod{R}$, so $(\F{Q}, \Rep{Q}{\Cot{R}})$ is a complete hereditary cotorsion pair in $\Rep{Q}{R}$ by Lemmas \ref{induced cp} and \ref{flat}. Thus one has $\Cot{Q}=\Rep{Q}{\Cot{R}}$.
\end{prf*}

\section{A definition of tensor products of representations}\label{tp}
\noindent
In this section, inspired by the work of Salarian and Vahed \cite{SV2016}, we give a definition of tensor product functors in the category of representations and collect some basic results, which are used in the rest of the paper.

\begin{bfhpg}[\bf Tensor products of representations]\label{tproduct}
Let $X\in\Rep{Q}{R}$ be a representation of quivers.
We construct for any $\ZZ$-module $G$, a representation $\Hom{X}{G}\in\Rep{\Qop}{\Rop}$ as follows.
\begin{itemlist}
\item For each vertex $i\in \Qop_0$, set $\Hom{X}{G}(i)=\Hom[\ZZ]{X(i)}{G}$;
\item For each arrow $a: i\to j$ in $\Qop_1$,
      define $\Hom{X}{G}(a)=\Hom[\ZZ]{X(\aop)}{G}: \Hom{X}{G}(i)\to \Hom{X}{G}(j)$,
      where $\aop: j \to i$ is an arrow in $Q$.
\end{itemlist}

It is evident that $\Hom{X}{-}$ is a functor from the category $\Mod{\ZZ}$ of all $\ZZ$-modules to $\Rep{\Qop}{\Rop}$.
This functor is left exact and preserves arbitrary products, so it has a left adjoint from $\Rep{\Qop}{\Rop}$ to $\Mod{\ZZ}$,
which is denoted $\tp{-}{X}$ and will play the role of the tensor product in our work.
\end{bfhpg}

The next result is clear from the definition, which is used frequently in this paper.

\begin{thm}\label{adjoint}
Let $X$ be in $\Rep{Q}{R}$ and $Y$ in $\Rep{\Qop}{\Rop}$. Then for each $G\in\Mod{\ZZ}$ there is a natural isomorphism
\begin{equation*}
\Hom[\ZZ]{\tp{Y}{X}}{G}\is\Hom[\Qop]{Y}{\Hom{X}{G}}\:.
\end{equation*}
\end{thm}

In the following, we collect some basic results on tensor products. We always let $X^+$ denote $\Hom{X}{\mathbb{Q}/\ZZ}\in\Rep{\Qop}{\Rop}$.

\begin{prp}
Let $X$ be in $\Rep{Q}{R}$ and $Y$ in $\Rep{\Qop}{\Rop}$. Then there is a natural isomorphism
\begin{equation*}
\tp{Y}{(\oplus_{p}X_p)}\is \oplus_{p}(\tp{Y}{X_p})\:.
\end{equation*}
\end{prp}
\begin{prf*}
By Theorem \ref{adjoint} one has
  \begin{align*}
    (\tp{Y}{(\oplus_{p}X_p)})^+
    &\dis \Hom[\Qop]{Y}{(\oplus_{p}X_p)^+}\\
    &\dis \Pi_{p}\Hom[\Qop]{Y}{(X_p)^+}\\
    &\dis \Pi_{p}(\tp[Q]{Y}{X_p})^+\\
    &\dis (\oplus_{p}(\tp[Q]{Y}{X_p}))^+
  \end{align*}
Thus there is an isomorphism $\tp{Y}{(\oplus_{p}X_p)}\is \oplus_{p}(\tp{Y}{X_p})$.
\end{prf*}

\begin{lem}\label{swap}
Let $X$ be in $\Rep{Q}{R}$ and $Y$ in $\Rep{\Qop}{\Rop}$. Then for each $G\in\Mod{\ZZ}$ there is a natural isomorphism
\begin{equation*}
\Hom[\Qop]{Y}{\Hom{X}{G}}\is\Hom[Q]{X}{\Hom{Y}{G}}\:.
\end{equation*}
\end{lem}
\begin{prf*}
Let $f=\{f(i):Y(i)\to\Hom{X}{G}(i)\}_{i\in Q_0}$ be in $\Hom[\Qop]{Y}{\Hom{X}{G}}$. For each $i\in Q_0$, the swap map
$$\zeta^{Y(i)GX(i)}: \Hom[\Rop]{Y(i)}{\Hom[\ZZ]{X(i)}{G}}\to \Hom[R]{X(i)}{\Hom[\ZZ]{Y(i)}{G}}$$
is an isomorphism with the inverse $\zeta^{X(i)GY(i)}$. So $\{\zeta^{Y(i)GX(i)}(f(i))\}_{i\in Q_0}$ is a morphism from $Y$ to $\Hom{X}{G}$. Set
$$\zeta^{YGX}: \Hom[\Qop]{Y}{\Hom{X}{G}}\to\Hom[Q]{X}{\Hom{Y}{G}}$$
with $\zeta^{YGX}(f)=\{\zeta^{Y(i)GX(i)}(f(i))\}_{i\in Q_0}$. Then it has an inverse $\zeta^{XGY}$. Thus one gets the isomorphism in the statement.
\end{prf*}

\begin{prp}
Let $X$ be in $\Rep{Q}{R}$ and $Y$ in $\Rep{\Qop}{\Rop}$. Then there is a natural isomorphism
\begin{equation*}
\tp{Y}{X}\is \tp[\Qop]{X}{Y}\:.
\end{equation*}
\end{prp}
\begin{prf*}
By Propositions \ref{adjoint} and Lemma \ref{swap} one has
  \begin{align*}
    (\tp{Y}{X})^+
    &\dis \Hom[\Qop]{Y}{X^+}\\
    &\dis \Hom[Q]{X}{Y^+}\\
    &\dis (\tp[\Qop]{X}{Y})^+\:.
  \end{align*}
So $\tp{Y}{X}$ is isomorphic to $\tp[\Qop]{X}{Y}$.
\end{prf*}

The next lemma is from Enochs, Estrada and Garc\'{i}a Rozas \corcite[6.7]{EEGR09}.

\begin{lem}\label{flatinj}
Let $Q$ be a left rooted quiver and $X\in\Rep{Q}{R}$. Then $X$ is flat if and only if $X^+$ is injective in $\Rep{\Qop}{\Rop}$.
\end{lem}

We end this section with the following result,
which may assert that our definition of tensor products of representations is reasonable.

\begin{thm}
Let $Q$ be a left rooted quiver and $X\in\Rep{Q}{R}$. Then $X$ is flat if and only if the functor $\tp{-}{X}$ is exact.
\end{thm}
\begin{prf*}
We assume that $X$ is flat and fix an exact sequence $0 \to Y' \to Y \to Y'' \to 0$ in $\Rep{\Qop}{\Rop}$. By Lemma \ref{flatinj}, $X^+$ is injective in $\Rep{\Qop}{\Rop}$, so the sequence
$$0 \to \Hom[\Qop]{Y''}{X^+} \to \Hom[\Qop]{Y}{X^+} \to \Hom[\Qop]{Y'}{X^+} \to 0$$
is exact. It follows from Theorem \ref{adjoint} that the sequence
$$0 \to \tp{Y'}{X} \to \tp{Y}{X} \to \tp{Y''}{X} \to 0.$$
So the functor $\tp{-}{X}$ is exact. Conversely, it is easy to see that $X^+$ is injective, so $X$ is flat by Lemma \ref{flatinj}.
\end{prf*}

\section{Gorenstein flat representations of quivers}\label{Gfd}
\noindent
In this section we prove Theorem A in the introduction. We begin with the following definitions.

\begin{bfhpg}[\bf Gorenstein injective representations of quivers]
A representation $Y\in\Rep{\Qop}{\Rop}$ is called \emph{Gorenstein injective} if there is an exact sequence
$$\cdots \to I_1 \to I_0 \to I_{-1} \to \cdots$$
in $\Rep{\Qop}{\Rop}$ with each $I_i$ injective for all $i\in\ZZ$, such that $X\is\Coker{(I_1 \to I_0)}$ and the sequence remains exact after applying the functor $\Hom[\Qop]{E}{-}$ for every injective representation $E\in\Rep{\Qop}{\Rop}$.

It follows from Eshraghi, Hafezi and Salarian \thmcite[3.5.1]{EHS13} that if $Q$ is a left rooted quiver (in this case $\Qop$ is a right rooted quiver) then $Y\in\Rep{\Qop}{\Rop}$ is Gorenstein injective if and only if $Y$ is in $\Psi(\GI{\Rop})$, that is, for each $i\in Q_{0}$ the homomorphism $\psi_i^Y$ is surjective, and the $\Rop$-modules $Y(i)$ and $\Ker{\psi_i^Y}$ are Gorenstein injective.
\end{bfhpg}

\begin{dfn}\label{gfqr}
A representation $X\in\Rep{Q}{R}$ is called \emph{Gorenstein flat} if there is an exact sequence
$$\cdots \to F_1 \to F_0 \to F_{-1} \to \cdots$$
in $\Rep{Q}{R}$ with each $F_i$ flat for all $i\in\ZZ$, such that $X\is\Coker{(F_1 \to F_0)}$ and the sequence remains exact after applying the functor $\tp{E}{-}$ for every injective representation $E\in\Rep{\Qop}{\Rop}$.
The subcategory of all Gorenstein flat representations in $\Rep{Q}{R}$ is denoted $\GF{Q}$.
\end{dfn}

\begin{lem}\label{onlyif}
Let $Q$ be a left rooted quiver and $X\in\Rep{Q}{R}$. If $X$ is Gorenstein flat, then $X$ is in $\Phi(\GF{R})$, that is, for each vertex $i\in Q_0$ the homomorphism $\varphi_i^X: \oplus_{a\in Q_{1}^{\ast\to i}}X(s(a))\to X(i)$ is injective, and the $R$-modules $X(i)$ and $\C_i(X)$ are Gorenstein flat.
\end{lem}
\begin{prf*}
Since $X$ is Gorenstein flat, there exists an exact sequence
$$\bF=\cdots \to F_1 \to F_0 \to F_{-1} \to \cdots$$
in $\Rep{Q}{R}$ with each $F_j$ flat for all $j\in\ZZ$, such that $X\is\Coker{(F_1 \to F_0)}$ and the sequence $\tp{E}{\bF}$ is exact for every injective representation $E\in\Rep{\Qop}{\Rop}$. Fix an integer $j$. By Enochs, Oyonarte and Torrecillas \prpcite[3.4]{EOT04}, for each $i\in Q_0$ the homomorphism $\varphi_i^{F_j}: \oplus_{a\in Q_{1}^{\ast\to i}}F_j(s(a))\to F_j(i)$ is injective, and the $R$-modules $F_j(i)$ and $\C_i(F_j)$ are flat. Consider the following commutative diagram with exact rows and columns:
$$\begin{array}{lcr}
\xymatrix@C=15pt@R=20pt{
  & \vdots \ar[d] & \vdots \ar[d] & \vdots \ar[d] \\
  0 \ar[r] & \oplus_{a\in Q_{1}^{\ast\to i}}F_1(s(a)) \ar[d]\ar[r] & F_1(i) \ar[d]\ar[r] & \C_i(F_1) \ar[d]\ar[r] & 0 \\
  0 \ar[r] & \oplus_{a\in Q_{1}^{\ast\to i}}F_0(s(a)) \ar[d]\ar[r] & F_0(i) \ar[d]\ar[r] & \C_i(F_0) \ar[d]\ar[r] & 0 \\
  0 \ar[r] & \oplus_{a\in Q_{1}^{\ast\to i}}F_{-1}(s(a)) \ar[d]\ar[r] & F_{-1}(i) \ar[d]\ar[r] & \C_i(F_{-1}) \ar[d]\ar[r] & 0. \\
  & \vdots & \vdots & \vdots  }
  \end{array} \eqno(\ast)$$
So the induced sequence $0 \to \oplus_{a\in Q_{1}^{\ast\to i}}X(s(a)) \xra{\varphi_i^X} X(i) \to \C_i(X) \to 0$ is exact. For every injective representation $E\in\Rep{\Qop}{\Rop}$, one has $\Hom[\Qop]{E}{\bF^+}\is(\tp{E}{\bF})^+$ is exact by Proposition \ref{adjoint}, so \thmcite[3.1.6(i)]{EHS13} yields that $\bF(i)^+$ is a totally acyclic complex of injective $\Rop$-modules; as $\Qop$ is right rooted. Thus for each injective $\Rop$-module $I$, the complex $\tp[R]{I}{\bF(i)}$ is acyclic. Hence $X(i)$ is Gorenstein flat. The complex $\tp[R]{I}{(\oplus_{a\in Q_{1}^{\ast\to i}}\bF(s(a)))}$ is acyclic and each $\C_i(F_j)$ is flat for $j\in\ZZ$, so from the diagram ($\ast$) one gets that the complex $\tp[R]{I}{\C_i(\bF)}$ is acyclic. Thus $\C_i(X)$ is Gorenstein flat.
\end{prf*}

The next result paves the way for Lemma \ref{cogenerator} which is used in the proof of Theorem A in the introduction;
it is proved using the Snake Lemma.

\begin{lem}\label{diagram}
Let $\sf W$ and $\sf X$ be subcategories of $\Mod{R}$, and let $0 \to X_1 \to X_2 \to X_3 \to 0$ be an exact sequence in $\sf X$. If there exist two exact sequence $0 \to X_1 \to W_1 \to X_1' \to 0$ and $0 \to X_3 \to W_3 \to X_3' \to 0$ with $W_1,W_3\in\sf W$ and $X_1', X_3' \in \sf X$ such that $\Ext{1}{X_3}{W_1}=0$, then there exists a commutative diagram
$$\xymatrix@C=15pt@R=15pt{
  & 0 \ar[d] & 0 \ar[d] & 0 \ar[d] \\
  0 \ar[r] & X_1 \ar[d]\ar[r] & X_2 \ar[d]\ar[r] & X_3 \ar[d]\ar[r] & 0 \\
  0 \ar[r] & W_1 \ar[d]\ar[r] & W_1\oplus W_3 \ar[d]\ar[r] & W_3 \ar[d]\ar[r] & 0 \\
  0 \ar[r] & X_1' \ar[d]\ar[r] & X_2' \ar[d]\ar[r] & X_3' \ar[d]\ar[r] & 0 \\
  & 0 & 0 & 0  }$$
  with exact rows and columns.
\end{lem}

\begin{ipg}
Let $\sf W$ and $\sf X$ be subcategories of an abelian category. Recall from Auslander and Buchweitz \cite{MAsROB89} that $\sf W$ is a \emph{cogenerator} for $\sf X$ if $\sf W\subseteq\sf X$, and for each $X\in\sf X$ there exists an exact sequence $0 \to X \to W \to X' \to 0$ with $W\in\sf W$ and $X'\in\sf X$.
\end{ipg}

The following lemma is essentially based on the construction given by Odaba\c{s}{\i} \lemcite[4.5]{Od19};
one refers to \ref{pf1} for a proof over a more simple quiver with 4 vertices.

\begin{lem}\label{cogenerator}
Assume that $\sf W$ and $\sf X$ are subcategories of $\Mod{R}$ satisfying
\begin{itemlist}
\item $\sf X$ is closed under extensions;
\item $(\sf W, {\sf W}^{\perp})$ is a complete cotorsion pair;
\item $\Ext{1}{X}{W}=0$ for all $X\in\sf X$ and $W\in {\sf W\cap{\sf W}^{\perp}}$.
\end{itemlist}
Let $Q$ be a left rooted quiver. If\: $\sf W$ is a cogenerator for $\sf X$, then $\Phi(\sf W)$ is a cogenerator for $\Phi(\sf X)$.
\end{lem}
\begin{prf*}
Let $\{V_{\alpha}\}$ be the transfinite sequence of subsets of $Q_0$. Since $Q$ is left rooted, one has $Q_0=V_{\lambda}$ for some ordinal $\lambda$. Let $X$ be in $\Phi(\sf X)$. For every ordinal $\alpha\leq\lambda$, we define a representation $X_{\alpha}$ as follows:
\begin{equation*}
  X_{\alpha}(i)=
  \begin{cases}
    X(i) & \text{if $i\in V_{\alpha}$}\,,\\
    0 & \text{if $i\notin V_{\alpha}$}\,.
  \end{cases}
\end{equation*}
For an arrow $a: j\to k$ in $Q$, the morphism $X_\alpha(a)$ is defined as
\begin{equation*}
  X_{\alpha}(a)=
  \begin{cases}
    X(a) & \text{if $j,k\in V_{\alpha}$}\,,\\
    0 & \text{otherwise}\,.
  \end{cases}
\end{equation*}
It is clear that $X_{\lambda}=X$. Next for each ordinal $\alpha\leq\lambda$, we construct an exact sequence
\begin{equation*}
\bE_{\alpha}:\ 0 \to X_{\alpha} \to W_{\alpha} \to Y_{\alpha} \to 0
\end{equation*}
in $\Rep{Q}{R}$ such that the following conditions hold:
\begin{prt}
\item For each $i\notin V_{\alpha}$, $W_\alpha(i)=0$, and for all $i\in V_{\alpha}$ the homomorphism $\varphi_i^{W_{\alpha}}$ is injective with both $W_\alpha(i)$ and $\C_i(W_{\alpha})$ belong to $\sf W$;
\item For each $i\notin V_{\alpha}$, $Y_\alpha(i)=0$, and for all $i\in V_{\alpha}$ the homomorphism $\varphi_i^{Y_{\alpha}}$ is injective with both $Y_\alpha(i)$ and $\C_i(Y_{\alpha})$ belong to $\sf X$;
\item If $\alpha<\alpha'\leq\lambda$, then there is a morphism from $\bE_{\alpha'}$ to $\bE_{\alpha}$, and $\bE_{\alpha'}(i)=\bE_{\alpha}(i)$ for all $i\notin V_{\alpha'}\setminus V_{\alpha}$.
\end{prt}
We do this by transfinite induction. The case where $\alpha=0$ holds clearly as $V_0=\emptyset$.
For each $i\in V_1$, there exists an exact sequence $0 \to X(i) \to W^i \to Y^i \to 0$ with $W^i\in\sf W$ and $X^i\in\sf X$, as $\sf W$ is a cogenerator for $\sf X$. We define representations $W_{1}$ and $Y_1$ as follows:
\begin{equation*}
  W_{1}(i)=
  \begin{cases}
    W^i & \text{if $i\in V_{1}$}\\
    0 & \text{if $i\notin V_{1}$}
  \end{cases}
  \ \mathrm{and}\
    Y_{1}(i)=
  \begin{cases}
    Y^i & \text{if $i\in V_{1}$}\\
    0 & \text{if $i\notin V_{1}$}\,.
  \end{cases}
\end{equation*}
For an arrow $a: j\to k$ in $Q$, we let $W_1(a)$ and $Y_1(a)$ be zero. Then there is an exact sequence $\bE_1: 0 \to X_1 \to W_1 \to Y_1 \to 0$ in $\Rep{Q}{R}$ satisfying the desired conditions, as there is no arrow $a$ with $t(a)\in V_1$.

Suppose that $\alpha+1$ is a successor ordinal and we have $\bE_{\alpha}$. Next we construct $\bE_{\alpha+1}$.

Let $i\in V_{\alpha+1}\setminus V_{\alpha}$. Since $X$ is in $\Phi(\sf X)$, the sequence
$$0 \to \oplus_{a\in Q_{1}^{\ast\to i}}X(s(a))\xra{\varphi_i^{X}} X(i) \to \C_i(X) \to 0$$
is exact. For all $a\in Q_{1}^{\ast\to i}$, one has $s(a)\in V_{\alpha}$ (see \ref{rooted}), so $X(s(a))=X_\alpha(s(a))$. Thus there is an exact sequence
$$\ 0 \to \oplus_{a\in Q_{1}^{\ast\to i}}X(s(a)) \to \oplus_{a\in Q_{1}^{\ast\to i}}W_{\alpha}(s(a)) \to \oplus_{a\in Q_{1}^{\ast\to i}}Y_{\alpha}(s(a)) \to 0.$$
By the assumption, $(\sf W, {\sf W}^{\perp})$ is a complete cotorsion pair, so there is an exact sequence $0 \to \oplus_{a\in Q_{1}^{\ast\to i}}W_{\alpha}(s(a)) \to U^i \to C' \to 0$ with $U^i\in{\sf W}^{\perp}$ and $C'\in\sf W$. Since $W_{\alpha}(s(a))$ is in $\sf W$, the module $\oplus_{a\in Q_{1}^{\ast\to i}}W_{\alpha}(s(a))$ is so. Hence $U^i$ is in ${\sf W}^{\perp}\cap \sf W$. Consider the following push-out diagram:
$$\xymatrix@C=10pt@R=20pt{&&0\ar[d] &0\ar[d]\\
  0 \ar[r] & \oplus_{a\in Q_{1}^{\ast\to i}}X(s(a)) \ar[r]^{}\ar@{=}[d] & \oplus_{a\in Q_{1}^{\ast\to i}}W_{\alpha}(s(a)) \ar[r]^{}\ar[d] & \oplus_{a\in Q_{1}^{\ast\to i}}Y_{\alpha}(s(a)) \ar[r]\ar[d] & 0\\
  0 \ar[r] & \oplus_{a\in Q_{1}^{\ast\to i}}X(s(a)) \ar[r]^{} & U^i \ar[r]^{}\ar[d] & T^i \ar[r]\ar[d] & 0.\\
  && C' \ar@{=}[r]\ar[d] & C' \ar[d]\\
  && 0 & 0 }$$
Since $\sf W$ is a cogenerator for $\sf X$ and $\C_i(X)\in\sf X$, one has an exact sequence
$$0 \to \C_i(X) \to S^i \to Z^i \to 0$$
with $S^i \in \sf W$ and $Z^i\in\sf X$. By the assumption one has $\Ext{1}{\C_i(X)}{U^i}=0$. Thus by Lemma \ref{diagram}, one gets a commutative diagram with exact rows and columns
$$\xymatrix@C=15pt@R=15pt{
  & 0 \ar[d] & 0 \ar[d] & 0 \ar[d] \\
  0 \ar[r] & \oplus_{a\in Q_{1}^{\ast\to i}}X(s(a)) \ar[d]\ar[r]^{\ \ \ \ \ \ \ \ \ \varphi_i^{X}} & X(i) \ar[d]\ar[r] & \C_i(X) \ar[d]\ar[r] & 0 \\
  0 \ar[r] & U^i \ar[d]\ar[r] & W^i \ar[d]\ar[r] & S^i \ar[d]\ar[r] & 0 \\
  0 \ar[r] & T^i \ar[d]\ar[r] & Y^i \ar[d]\ar[r] &Z^i \ar[d]\ar[r] & 0, \\
  & 0 & 0 & 0  }$$
where $W^i=U^i\oplus S^i\in\sf W$ and $Y^i\in\sf X$, as $\sf W$ and $\sf X$ are closed under extensions.

We define a representation $W_{\alpha+1}$ as follows:
\begin{equation*}
  W_{\alpha+1}(i)=
  \begin{cases}
    W^i & \text{if $i\in V_{\alpha+1}\backslash V_{\alpha}$}\,,\\
    W_\alpha(i) & \text{if $i\in V_{\alpha}$}\,,\\
    0 & \text{if $i\notin V_{\alpha+1}$}\,.
  \end{cases}
\end{equation*}
For an arrow $a: j\to k$ in $Q$, the morphism $W_{\alpha+1}(a)$ is defined as follows:
\begin{rqm}
\item[--] If $k\in V_{\alpha+1}\backslash V_{\alpha}$, then $j\in V_{\alpha}$, and so $W_{\alpha+1}(j)=W_{\alpha}(j)$ and $W_{\alpha+1}(k)=W^{k}$. In this case, $W_{\alpha+1}(a)$ is the compositions
    $$W_{\alpha}(j)\hookrightarrow\oplus_{a\in Q_{1}^{\ast\to k}}W_{\alpha}(s(a))\hookrightarrow U^k\hookrightarrow W^k.$$
\item[--] If $k\in V_{\alpha}$, then $j\in V_{\alpha-1}\subseteq V_{\alpha}$, and so $W_{\alpha+1}(j)=W_{\alpha}(j)$ and $W_{\alpha+1}(k)=W_{\alpha}(k)$. In this case, $W_{\alpha+1}(a)=W_{\alpha}(a)$.
\item[--] If $k\notin V_{\alpha+1}$, then $W_{\alpha+1}(k)=0$. In this case, $W_{\alpha+1}(a)=0$.
\end{rqm}
It is obvious that $W_{\alpha+1}(i)$ is in $\sf W$ for all $i\in V_{\alpha+1}$.

Similarly, we define a representation $Y_{\alpha+1}$ as follows:
\begin{equation*}
  Y_{\alpha+1}(i)=
  \begin{cases}
    Y^i & \text{if $i\in V_{\alpha+1}\backslash V_{\alpha}$}\,,\\
    Y_\alpha(i) & \text{if $i\in V_{\alpha}$}\,,\\
    0 & \text{if $i\notin V_{\alpha+1}$}\,.
  \end{cases}
\end{equation*}
For an arrow $a: j\to k$ in $Q$, the morphism $Y_{\alpha+1}(a)$ is defined as follows:
\begin{rqm}
\item[--] If $k\in V_{\alpha+1}\backslash V_{\alpha}$, then $j\in V_{\alpha}$, and so $Y_{\alpha+1}(j)=Y_{\alpha}(j)$ and $Y_{\alpha+1}(k)=Y^{k}$. In this case, $Y_{\alpha+1}(a)$ is the compositions
    $$Y_{\alpha}(j)\hookrightarrow\oplus_{a\in Q_{1}^{\ast\to k}}Y_{\alpha}(s(a))\hookrightarrow T^k\hookrightarrow Y^k.$$
\item[--] If $k\in V_{\alpha}$, then $j\in V_{\alpha-1}\subseteq V_{\alpha}$, and so $Y_{\alpha+1}(j)=Y_{\alpha}(j)$ and $Y_{\alpha+1}(k)=Y_{\alpha}(k)$. In this case, $Y_{\alpha+1}(a)=Y_{\alpha}(a)$.
\item[--] If $k\notin V_{\alpha+1}$, then $Y_{\alpha+1}(k)=0$. In this case, $Y_{\alpha+1}(a)=0$.
\end{rqm}
It is obvious that $Y_{\alpha+1}(i)$ is in $\sf X$  for all $i\in V_{\alpha+1}$.

Thus there exists an exact sequence
\begin{equation*}
\bE_{\alpha+1}:\ 0 \to X_{\alpha+1} \to W_{\alpha+1} \to Y_{\alpha+1} \to 0
\end{equation*}
in $\Rep{Q}{R}$; the corresponding commutative diagram for an arrow $a$ in $Q$ can be checked by the above construction. It is easy to see that there is a morphism from $\bE_{\alpha+1}$ to $\bE_{\alpha}$, and $\bE_{\alpha+1}(i)=\bE_{\alpha}(i)$ for all $i\notin V_{\alpha+1}\setminus V_{\alpha}$. In the following, we prove that for all $i\in V_{\alpha+1}$, the homomorphism $\varphi_i^{W_{\alpha+1}}$ is injective with $\C_i(W_{\alpha+1})\in\sf W$, and the homomorphism $\varphi_i^{Y_{\alpha+1}}$ is injective with $\C_i(Y_{\alpha+1})\in\sf X$. We only prove the case for $W_{\alpha+1}$.

Let $i\in V_{\alpha+1}$. We deal with the following two situations:
\begin{rqm}
\item If $i\in V_{\alpha+1}\setminus V_{\alpha}$, then $s(a)\in V_{\alpha}$ for all $a\in Q_{1}^{\ast\to i}$, and so $W_{\alpha+1}(s(a))=W_{\alpha}(s(a))$ and $W_{\alpha+1}(i)=W^{i}$. In this case, from the construction of $W_{\alpha+1}$, one sees that $\varphi_i^{W_{\alpha+1}}: \oplus_{a\in Q_{1}^{\ast\to i}}W_{\alpha+1}(s(a))\to W_{\alpha+1}(i)$ is the compositions
    $\oplus_{a\in Q_{1}^{\ast\to i}}W_{\alpha}(s(a))\hookrightarrow U^i\hookrightarrow W^i$, which is a monomophism. Using the Snake Lemma, one gets the next commutative diagram with exact rows and columns:
    $$\xymatrix@C=20pt@R=20pt{&&0\ar[d] &0\ar[d]\\
  0 \ar[r] & \oplus_{a\in Q_{1}^{\ast\to i}}W_{\alpha}(s(a)) \ar[r]^{}\ar@{=}[d] & U^i \ar[r]\ar[d] & C' \ar[r]\ar[d] & 0\\
  0 \ar[r] & \oplus_{a\in Q_{1}^{\ast\to i}}W_{\alpha}(s(a)) \ar[r]^{\ \ \ \ \ \ \ \ \ \ \varphi_i^{W_{\alpha+1}}} & W^i \ar[r]\ar[d] & \C_i(W_{\alpha+1}) \ar[r]\ar[d] & 0.\\
  && S^i \ar@{=}[r]\ar[d] & S^i \ar[d]\\
  && 0 & 0 }$$
  Since $C'$ and $S^i$ are in $\sf W$, so is $\C_i(W_{\alpha+1})$ as $\sf W$ is closed under extensions.
\item If $i\in V_{\alpha}$, then $s(a)\in V_{\alpha-1}\subseteq V_{\alpha}$ for all $a\in Q_{1}^{\ast\to i}$, and so $W_{\alpha+1}(s(a))=W_{\alpha}(s(a))$ and $W_{\alpha+1}(i)=W_{\alpha}(i)$. In this case, from the construction of $W_{\alpha+1}$, one sees that $\varphi_i^{W_{\alpha+1}}=\varphi_i^{W_{\alpha}}$ is a monomorphism. Hence $\C_i(W_{\alpha+1})=\C_i(W_{\alpha})$ is in $\sf W$.
\end{rqm}

Finally, if $\alpha\leq\lambda$ is a limit ordinal and if $\bE_\beta$ is constructed for all $\beta<\alpha$. Next we construct $\bE_\alpha$. In this case one has $V_\alpha=\cup_{\beta<\alpha}V_{\beta}$. If $i\in V_\alpha$, then $i\in V_\beta$ for some ordinal $\beta<\alpha$, and so for all ordinal $\alpha>\beta'\geq\beta$ one has $\bE_{\beta'}(i)=\bE_{\beta}(i)$ as $i\notin V_{\beta'}\backslash V_{\beta}$. If $i\notin V_\alpha$, then $i\notin V_\beta$ for each $\beta<\alpha$, and so $\bE_{\beta}(i)$ is an exact sequence of zero representations. We let $\bE_\alpha=\lim_{\beta<\alpha}\bE_\beta$. Then one has
\begin{equation*}
  \bE_{\alpha}(i)=
  \begin{cases}
    \bE_{\beta}(i)\ \text{for some}\ \beta<\alpha  & \text{if $i\in V_{\alpha}$}\,,\\
    0 & \text{if $i\notin V_{\alpha}$}\,.
  \end{cases}
\end{equation*}
Thus $\bE_{\alpha}$ is an exact sequence in $\Rep{Q}{R}$ satisfying the desired conditions.

As a consequence one gets an exact sequence $\ 0 \to X \to W \to Y \to 0$ in $\Rep{Q}{R}$ with $W\in\Phi(\sf W)$ and $Y\in\Phi(\sf X)$. So $\Phi(\sf W)$ is a cogenerator for $\Phi(\sf X)$.
\end{prf*}

\begin{lem}\label{fcogenerator}
Let $Q$ be a left rooted quiver. Then $\Phi(\F{R})$ is a cogenerator for $\Phi(\GF{R})$.
\end{lem}
\begin{prf*}
It is known that $\F{R}$ is a cogenerator for $\GF{R}$ and $(\F{R},\F{R}^{\perp})$ is a complete cotorsion pair; see Bican, El Bashir and Enochs \prpcite[2]{BEE-01}. It follows from \corcite[4.12]{SS20} that $\GF{R}$ is closed under extensions and $\Ext{>0}{G}{W}=0$ for each $G\in\GF{R}$ and each $W\in\F{R}\cap\F{R}^{\perp}$. So by Lemma \ref{cogenerator}, $\Phi(\F{R})$ is a cogenerator for $\Phi(\GF{R})$.
\end{prf*}

\begin{lem}\label{character}
If $X$ is in $\Phi(\GF{R})$, then $X^+$ is in $\Psi(\GI{R})$.
\end{lem}
\begin{prf*}
For all $i\in Q_0$, one has $\Hom[\ZZ]{\varphi_{i}^{X}}{\mathbb{Q}/\ZZ}=\psi_i^{X^+}$, $\Hom[\ZZ]{X(i)}{\mathbb{Q}/\ZZ}=X^+(i)$ and $\Hom[\ZZ]{\C_i(X)}{\mathbb{Q}/\ZZ}=\mathrm{K}_i(X^+)$. It is known that if $M$ is a Gorenstein flat $R$-module then the $\Rop$-module $\Hom[\ZZ]{M}{\mathbb{Q}/\ZZ}$ is Gorenstein injective; see Holm \thmcite[3.6]{HHl04a}. Thus if $X$ is in $\Phi(\GF{R})$, then $X^+$ is in $\Psi(\GI{R})$.
\end{prf*}

\begin{bfhpg}[\bf Proof of Theorem A]
The ``only if" part holds by Lemma \ref{onlyif}. For the ``if" part we assume that for each $i\in Q_0$ the homomorphism $\varphi_i^X: \oplus_{a\in Q_{1}^{\ast\to i}}X(s(a))\to X(i)$ is injective, and the $R$-modules $X(i)$ and $\C_i(X)$ are Gorenstein flat, that is, $X$ is in $\Phi(\GF{R})$.

By Lemma \ref{fcogenerator}, there is an exact sequence $0 \to X \to F_{-1} \to X_{-1} \to 0$ in $\Rep{Q}{R}$ with $F_{-1}\in\Phi(\F{R})$ and $X_{-1}\in\Phi(\GF{R})$. So the sequence $0 \to X_{-1}^+ \to F_{-1}^+ \to X^+ \to 0$ is exact. By Lemma \ref{character}, $X_{-1}^+$ is in $\Psi(\GI{R})$, so $X_{-1}^+$ is Gorenstein injective; see \thmcite[3.5.1(a)]{EHS13}. Hence for each injective representation $I\in\Rep{\Qop}{\Rop}$, the sequence
$$0 \to \Hom[\Qop]{I}{X_{-1}^+} \to \Hom[\Qop]{I}{F_{-1}^+} \to \Hom[\Qop]{I}{X^+} \to 0$$
is exact. It follows from Proposition \ref{adjoint} that the sequence
$$0 \to \tp{I}{X} \to \tp{I}{F_{-1}} \to \tp{I}{X_{-1}} \to 0$$
is exact. Since $X_{-1}$ is in $\Phi(\GF{R})$, continuing the above process one gets an exact sequence
\begin{equation*}
\tag{$\ast$}
0 \to X \to F_{-1} \to F_{-2} \to \cdots
\end{equation*}
in $\Rep{Q}{R}$ with each $F_i\in\Phi(\F{R})$ (so $F_i\in\F{Q}$ by \thmcite[3.7]{EOT04}) such that it remains exact after applying the functor $\tp{I}{-}$ for each injective representation $I\in\Rep{\Qop}{\Rop}$.

On the other hand, there is an exact sequence $0 \to X_1 \to F_0 \to X \to 0$ in $\Rep{Q}{R}$ with $F_0$ flat. By \thmcite[3.7]{EOT04}, $F_0$ is in $\Phi(\F{Q})$, so it follows from \corcite[4.12]{SS20} that $X_1$ is in $\Phi(\GF{Q})$. The same argument as above yields that the sequence $0 \to \tp{I}{X_{1}} \to \tp{I}{F_{0}} \to \tp{I}{X} \to 0$ is exact for each injective representation $I\in\Rep{\Qop}{\Rop}$. Continuing this process, one gets an exact sequence
\begin{equation*}
\tag{$\ast\ast$}
\cdots \to F_1 \to F_{0} \to X \to 0
\end{equation*}
in $\Rep{Q}{R}$ with each $F_i\in\F{Q}$, such that it remains exact after applying the functor $\tp{I}{-}$ for each injective representation $I\in\Rep{\Qop}{\Rop}$.

Assembling the sequences $(\ast)$ and $(\ast\ast)$, one gets an exact sequence
$$\cdots \to F_1 \to F_{0} \to F_{-1} \to F_{-2} \to \cdots$$
in $\Rep{Q}{R}$ with each $F_i$ flat for $i\in\ZZ$, such that $X\is\Coker{(F_1 \to F_0)}$ and the sequence remains exact after applying the functor $\tp{I}{-}$ for every injective representation $I\in\Rep{\Qop}{\Rop}$. This yields that $X$ is a Gorenstein flat representation in $\Rep{Q}{R}$.
\end{bfhpg}

\begin{cor}\label{gfpair}
Let $Q$ be a left rooted quiver. Then $\GF{Q}$ is closed under direct limits and $(\GF{Q}, \Rep{Q}{\GF{R}^\perp})$ is a perfect hereditary cotorsion pair in $\Rep{Q}{R}$ with $\GF{Q}\cap\Rep{Q}{\GF{R}^\perp}=\F{Q}\cap\Rep{Q}{\Cot{R}}.$
\end{cor}
\begin{prf*}
It follows from \corcite[4.12]{SS20} that $(\GF{R},\GF{R}^\perp)$ is a complete hereditary cotorsion pair in $\Mod{R}$. Thus by Lemma \ref{induced cp}  the pair
$(\Phi(\GF{R}), \Rep{Q}{\GF{R}^\perp})$ is a complete hereditary cotorsion pair in $\Rep{Q}{R}$. Since $\GF{R}$ is closed under direct limits again by \corcite[4.12]{SS20}, one gets that $\Phi(\GF{R})$ is so. Hence the cotorsion pair $(\Phi(\GF{R}), \Rep{Q}{\GF{R}^\perp})$ is perfect. Next we prove $\Phi(\GF{R})\cap\Rep{Q}{\GF{R}^\perp}=\Phi(\F{R})\cap\Rep{Q}{\Cot{R}}$.

According to \corcite[4.12]{SS20},
one has $\GF{R}\cap \GF{R}^\perp= \F{R}\cap\Cot{R}$.
This implies that the inclusion ``$\supseteq$" holds true. For the inclusion ``$\subseteq$", we let $X\in \Phi(\GF{R})\cap\Rep{Q}{\GF{R}^\perp}$.
Then for each $i\in Q_0$,
$X(i)\in \GF{R}\cap \GF{R}^\perp= \F{R}\cap\Cot{R}$.
To show $X\in \Phi(\F{R})\cap\Rep{Q}{\Cot{R}}$,
it suffices to prove that $\C_i(X)\in \F{R}$ for each $i\in Q_0$.
To this end, consider the following short exact sequence in $\Mod{R}$
\begin{center}
$0\to \oplus_{a\in Q_{1}^{\ast\to i}}X(s(a))\overset{\varphi_i^X}\longrightarrow X(i)\to \C_i(X)\to0$.
\end{center}
Note that $\oplus_{a\in Q_{1}^{\ast\to i}}X(s(a))$ is flat.
Then the flat dimension of the Gorenstein flat module $\C_i(X)$ is finite.
Hence one has $\C_i(X)\in \F{R}$.

Finally, by \thmcite[3.7]{EOT04} and Theorem A, one has $\Phi(\F{R})=\F{Q}$ and $\Phi(\GF{R})=\GF{Q}$. This finishes the proof.
\end{prf*}

We may use Theorem A to relate Gorenstein flat and Gorenstein injective representations.

\begin{cor}\label{gflatinj}
If $X$ is Gorenstein flat in $\Rep{Q}{R}$, then $X^+$ is Gorenstein injective in $\Rep{\Qop}{\Rop}$. The converse holds if $R$ is right coherent.
\end{cor}
\begin{prf*}
For all $i\in Q_0$, one has $\Hom[\ZZ]{\varphi_{i}^{X}}{\mathbb{Q}/\ZZ}=\psi_i^{X^+}$, $\Hom[\ZZ]{X(i)}{\mathbb{Q}/\ZZ}=X^+(i)$ and $\Hom[\ZZ]{\C_i(X)}{\mathbb{Q}/\ZZ}=\mathrm{K}_i(X^+)$. It is follows from \thmcite[3.6]{HHl04a} that if $M$ is a Gorenstein flat $R$-module then the $\Rop$-module $\Hom[\ZZ]{M}{\mathbb{Q}/\ZZ}$ is Gorenstein injective, and the converse holds if $R$ is right coherent. So using Theorem A and \thmcite[3.5.1]{EHS13} one can complete the proof.
\end{prf*}

\subsection*{\textbf{\em Projectively coresolved Gorenstein flat representations}} %%%%%%%%%%%%%%
In this subsection, we give a similar characterization for projectively coresolved Gorenstein flat representations of quivers, which is used in the proof of Corollary C in the introduction.

\begin{bfhpg}[\bf Projectively coresolved Gorenstein flat objects]
Recall from \cite{SS20} that an $R$-module $M$ is \emph{projectively coresolved Gorenstein flat} if there is an exact sequence
$$\cdots \to P^{-1} \to P^0 \to P^{1} \to \cdots$$ of projective $R$-modules with $M\is\Ker{(P^0 \to P^1)}$, such that it remains exact after applying the functor $\tp[R]{I}{-}$ for every injective $\Rop$-module $I$. The subcategory of all projectively coresolved Gorenstein flat $R$-modules is denoted $\PGF{R}$.

Similarly, a representation $X\in\Rep{Q}{R}$ is called \emph{projectively coresolved Gorenstein flat} if there is an exact sequence
$$\cdots \to P_1 \to P_0 \to P_{-1} \to \cdots$$
in $\Rep{Q}{R}$ with each $P_i$ projective for all $i\in\ZZ$, such that $X\is\Coker{(P_1 \to P_0)}$ and the sequence remains exact after applying the functor $\tp{E}{-}$ for every injective representation $E\in\Rep{\Qop}{\Rop}$. The subcategory of all projectively coresolved Gorenstein flat representations in $\Rep{Q}{R}$ is denoted $\PGF{Q}$.
\end{bfhpg}

We let $\Prj{Q}$ denote the subcategory of projective representations in $\Rep{Q}{R}$. Let $Q$ be a left rooted quiver.
Then by \thmcite[3.1]{EE05} one has $\Prj{Q}=\Phi(\Prj{R})$.
It follows from \thmcite[4.9]{SS20} that the pair $(\PGF{R}, \PGF{R}^{\perp})$ is a complete hereditary cotorsion pair
with $\PGF{R}\cap\PGF{R}^{\perp}=\Prj{R}$,
so by Lemma \ref{cogenerator} one gets that $\Prj{Q}=\Phi(\Prj{R})$ is a cogenerator for $\Phi(\PGF{R})$.
Thus the same argument as in Theorem A yields the next result.

\begin{thm}\label{PGF}
Let $Q$ be a left rooted quiver and $X\in\Rep{Q}{R}$.
Then $X$ is projectively coresolved Gorenstein flat if and only if $X$ is in $\Phi(\PGF{R})$,
that is, for each $i\in Q_0$ the homomorphism
$\varphi_i^X: \oplus_{a\in Q_{1}^{\ast\to i}}X(s(a))\to X(i)$ is injective,
and the $R$-modules $X(i)$ and $\C_i(X)$ are projectively coresolved Gorenstein flat.
\end{thm}

By \thmcite[4.4]{SS20}, each module in $\PGF{R}$ is Gorenstein projective. Therefore each module in $\PGF{R}$ that has finite projective dimension must be projective.
So using a similar proof as in Corollary \ref{gfpair}, one gets the next corollary.

\begin{cor}\label{pgfpair}
Let $Q$ be a left rooted quiver. Then $(\PGF{Q}, \Rep{Q}{\PGF{R}^\perp})$ is a complete hereditary cotorsion pair in $\Rep{Q}{R}$ with
$$\PGF{Q}\cap\Rep{Q}{\PGF{R}^\perp}=\Prj{Q}.$$
\end{cor}

\section{An application to model structures}\label{model}
\noindent
In this section we aim at proving Theorem B in the introduction.
The next result is a key to do that, which is from Gillespie \thmcite[1.1]{G15}.
Recall that a subcategory $\sf{W}$ of an abelian category is called \emph{thick}
provided that it is closed under direct summands, extensions, and taking kernels of epimorphisms and cokernels of
monomorphisms.

\begin{thm}\label{How to}
Let $\sf{A}$ be an abelian category.
Suppose that
$(\sf{C},\widetilde{\sf{F}})$ and $(\widetilde{\sf{C}},\sf{F})$
are two complete hereditary cotorsion pairs in $\mathsf{A}$ such that
$\widetilde{\mathsf{C}}\subseteq \mathsf{C}$ (or equivalently, $\widetilde{\sf{F}}\subseteq \sf{F})$ and $\mathsf{C}\cap\widetilde{\mathsf{F}}=\widetilde{\mathsf{C}}\cap\mathsf{F}$.
Then there exists an unique thick subcategory $\mathsf{W}$
for which $(\mathsf{C}, \mathsf{W}, \mathsf{F})$ forms a Hovey triple.
Moreover, this thick subcategory can be described as follows:
\begin{center}
$\mathsf{W}=\{M\in \mathsf{A}\,|$ there exists a short exact sequence
$0\to M\to A\to B\to 0$\\
in $\mathsf{A}$ with $A\in \widetilde{\mathsf{F}}$ and $B\in \widetilde{\mathsf{C}}$
$\}$.
\end{center}
\end{thm}

\begin{ipg}\label{induced CP}
Let $Q$ be a left rooted quiver and $\sf{A}$ an abelian category with enough injectives
which satisfies the axiom $\sf{AB4}$.
Suppose that
$(\sf{C},\W,\mathsf{F})$ is a hereditary Hovey triple on $\sf{A}$,
and the associated complete hereditary cotorsion pairs in $\sf{A}$ are
$(\sf{C},\widetilde{\sf{F}}=\W\cap\sf{F})$ and $(\widetilde{\sf{C}}=\sf{C}\cap\W,\sf{F})$.
Then by \thmcite[7.9]{HJ19} and \thmcite[4.6]{Od19},
we obtain two complete hereditary cotorsion pairs
\begin{center}
$(\Phi({\sf C}),\Rep{Q}{\widetilde{\sf{F}}})$ \quad and \quad $(\Phi(\widetilde{\mathsf{C}}),\Rep{Q}{\sf{F}})$
\end{center}
in $\Rep{Q}{\sf{A}}$.

It is evident that $\Phi(\widetilde{\mathsf{C}})\subseteq \Phi({\sf C})$ since $\widetilde{\mathsf{C}}\subseteq \mathsf{C}$.
Next we show that
\begin{center}
$\Phi({\sf C})\cap\Rep{Q}{\widetilde{\sf{F}}}=\Phi(\widetilde{\mathsf{C}})\cap\Rep{Q}{\sf{F}}$\footnote{This equality was proved by Dalezios in \cite[pp. 5]{Da20} under the assumption that $\W$ is closed under coproducts.}.
\end{center}
The nontrivial inclusion is
$\Phi({\sf C})\cap\Rep{Q}{\widetilde{\sf{F}}}\subseteq\Phi(\widetilde{\mathsf{C}})\cap\Rep{Q}{\sf{F}}$.
Let $X\in \Phi({\sf C})\cap\Rep{Q}{\widetilde{\sf{F}}}$.
Then $X(i)\in {\sf C}\cap \widetilde{\sf{F}}= \widetilde{\mathsf{C}}\cap\mathsf{F}$ for each $i\in Q_0$.
Hence it remains to show that $\C_i(X)$ in the exact sequence
$$0 \to \oplus_{a\in Q_{1}^{\ast\to i}}X(s(a))\xra{\varphi_i^{X}} X(i) \to \C_i(X) \to 0$$
belongs to $\widetilde{\mathsf{C}}$.
To this end, note that $\widetilde{\sf{F}}= {\sf W}\cap{\sf F}$.
It follows that $X(i)\in {\sf W}$ for each $i\in Q_0$.
On the other hand,
since ${\sf C}\cap \widetilde{\sf{F}}=\widetilde{\mathsf{C}}\cap\mathsf{F}$,
we see that $X(i)\in \widetilde{\mathsf{C}}$ as well for each $i\in Q_0$.
This yields $\oplus_{a\in Q_{1}^{\ast\to i}}X(s(a))\in \widetilde{\mathsf{C}}$.
In particular,
$\oplus_{a\in Q_{1}^{\ast\to i}}X(s(a))\in {\sf W}$
as $\widetilde{\mathsf{C}}= {\sf C}\cap {\sf W}$.
Since ${\sf W}$ is a thick subcategory of $\sf{A}$,
we conclude that $\C_i(X)\in {\sf W}$ as well.
Note that $\C_i(X)\in {\sf C}$.
It follows that $\C_i(X)\in {\sf C}\cap{\sf W}= \widetilde{\mathsf{C}}$, as desired.
\end{ipg}

According to what we discussed above,
we get the next result by Theorem \ref{How to}.

\begin{lem}\label{model structor}
Let $Q$ be a left rooted quiver and $\sf{A}$ an abelian category with enough injectives
which satisfies the axiom $\sf{AB4}$.
Suppose that
$(\sf{C},\W,\mathsf{F})$ is a hereditary Hovey triple on $\sf{A}$,
and the associated complete hereditary cotorsion pairs in $\sf{A}$ are
$(\sf{C},\widetilde{\sf{F}}=\W\cap\sf{F})$ and $(\widetilde{\sf{C}}=\sf{C}\cap\W,\sf{F})$.
Then there exists a hereditary Hovey triple
$(\Phi({\sf C}),\T,\Rep{Q}{\sf{F}})$
on $\Rep{Q}{\sf{A}}$ in which
\begin{center}
$\mathsf{T}=\{X\in \Rep{Q}{\sf{A}}\,|$ there exists a short exact sequence
$0\to X\to A\to B\to 0$\\
in \Rep{Q}{\sf{A}} with $A\in \Rep{Q}{\widetilde{\mathsf{F}}}$ and $B\in \Phi(\widetilde{{\sf C}})\}$.
\end{center}
\end{lem}

The following result,
which is from Gillespie \prpcite[2.4]{Gil162},
gives a clear characterization of the thick subcategory $\W$ in a Hovey triple $(\sf{C},\W,\mathsf{F})$
on an abelian category $\sf{A}$.

\begin{lem}\label{model structor111}
Let $\sf{A}$ be an abelian category.
Suppose that
$(\sf{C},\W,\mathsf{F})$ is a hereditary Hovey triple on $\sf{A}$.
Then the thick subcategory $\W$ is characterized in the next way:
\begin{center}
$\mathsf{W}=\{M\in \mathsf{A}\,|$ there exists a short exact sequence
$0\to M\to A\to B\to 0$\\
in $\mathsf{A}$ with $A\in \W\cap\sf{F}$ and $B\in \sf{C}\cap\W\}$.
\end{center}
\end{lem}

The next result shows that the subcategory of all trivial objects in the model structure
given in Proposition \ref{model structor}
coincides with the subcategory of $\Rep{Q}{\sf A}$
consisting of all $\mathsf{W}$-valued representations,
where $\W$ is described as in Lemma \ref{model structor111}.
One refers to \ref{pf2} for a proof over a more simple quiver with 4 vertices.

\begin{lem}\label{trivial objects}
Keep the assumptions and notation as in Lemma \ref{model structor}.
Then there is an equality
$$\mathsf{T}=\Rep{Q}{\W}.$$
\end{lem}

\begin{prf*}
The inclusion ``$\mathsf{T}\subseteq\Rep{Q}{\W}$" holds trivially.
Next we prove the inclusion ``$\Rep{Q}{\W}\subseteq\mathsf{T}$".
Let $X\in \Rep{Q}{\W}$.
Then for each $i\in Q_0$,
there exists a short exact sequence
$$0\to X(i) \overset{k(i)}\longrightarrow A(i) \overset{h(i)}\longrightarrow B(i) \to 0$$
in $\A$ with $A(i)\in \widetilde{\sf{F}}$ and $B(i)\in\widetilde{\sf{C}}$; see Lemma \ref{model structor111}.
Note that $\widetilde{\sf{F}}\subseteq \sf{F}$
and $(\widetilde{\sf{C}},\sf{F})$ is a cotorsion pair in $\A$, so
for each arrow $a: i\to j \in Q_1$, and each map
$$\xymatrix{X(i)\ar[rr]^{k(j)\circ X(a)}& & A(j)},$$ there exists a map $A(a):A(i)\to A(j)$ (and hence a map $B(a):B(i)\to B(j)$)
such that the following diagram with exact rows
$$\xymatrix{
0\ar[r] &X(i)\ar[r]^{k(i)}\ar[d]_{X(a)}& A(i) \ar@.[d]_{A(a)} \ar[r]^{h(i)}&  B(i) \ar@.[d]^{B(a)}\ar[r]&0  \\
0\ar[r] &X(j)\ar[r]^{k(j)}&  A(j)  \ar[r]^{h(j)} & B(j)\ar[r]&0}$$ is commutative.
Hence one gets an exact sequence
$$\bE:\,\,\,0\to X\to A\to B\to 0$$
in $\Rep{Q}{\A}$ such that $A\in \Rep{Q}{\widetilde{\sf{F}}}$ and $B\in \Rep{Q}{\widetilde{\sf{C}}}$.
But $B$ may fail in $\Phi(\widetilde{\sf{C}})$.
We will use transfinite induction to construct an exact sequence
$$\bE':\,\,\,0\to X\to A'\to B'\to 0$$
in $\Rep{Q}{\A}$ such that
$A'\in \Rep{Q}{\widetilde{\sf{F}}}$ and $B'\in\Phi(\widetilde{\sf{C}})$;
this yields that $X$ is in $\mathsf{T}$.

Let $\{V_{\alpha}\}$ be the transfinite sequence of subsets of $Q_0$. Since $Q$ is left rooted, one has $Q_0=V_{\lambda}$ for some ordinal $\lambda$. Next we construct a continuous direct $\lambda$-sequence
$$\{\,\bE_\alpha:\,\,\, 0\to X_\alpha \overset{k_\alpha}\longrightarrow A_\alpha\overset{h_\alpha}\longrightarrow B_\alpha\to0 \,\}_{\alpha\leqslant\lambda}$$
of short exact sequences in $\Rep{Q}{\A}$ satisfying the following conditions:

\begin{prt}
\item For every ordinal $0<\alpha\leqslant\lambda$, $X_\alpha=X$.
\item For every ordinal $0<\alpha\leq\lambda$ and $i\in Q_0\backslash V_\alpha$,
$A_\alpha(i) = A(i)\in \widetilde{\sf{F}}$ and $B_\alpha(i) = B(i)\in \widetilde{\sf{C}}$.
\item For every ordinal $\alpha\leqslant\lambda$ and $i\in V_\alpha$,
$A_\alpha(i)\in \widetilde{\sf{F}}$,  $\varphi_i^{B_\alpha}$ is a monomorphism,
$B_\alpha(i)\in \widetilde{\sf{C}}$ and $\C_i(B_\alpha)\in \widetilde{\sf{C}}$.
\item For every $0<\alpha<\alpha'\leqslant\lambda$,
$A_\alpha(i) = A_{\alpha'}(i)$ and $B_\alpha(i) = B_{\alpha'}(i)$
 for all $i\notin V_{\alpha'}\backslash V_\alpha$,
 and there exists the commutative diagram in $\Rep{Q}{\A}$
 $$\xymatrix{
\bE_\alpha:\quad    0\ar[r] & X_\alpha\ar[r]^{k_\alpha}\ar@{=}[d] &  A_\alpha \ar[d]_{f_{{\alpha'},\alpha}} \ar[r]^{h_\alpha}
&  B_\alpha \ar[d]_{{g_{{\alpha'},\alpha}}} \ar[r]&0  \\
\bE_{\alpha'}:\quad 0\ar[r] & X_{\alpha'}\ar[r]^{k_{\alpha'}}     &  A_{\alpha'}  \ar[r]^{h_{\alpha'}}
&  B_{\alpha'}\ar[r]                              &0}$$
such that both $f_{{\alpha'},\alpha}$ and $g_{{\alpha'},\alpha}$ are monomorphisms.
\end{prt}

Set $\bE_0= 0\to 0\to 0\to 0\to 0$ and $\bE_1 = \bE$.
Suppose that $\alpha+1$ is a successor ordinal and we have $\bE_{\alpha}= \,0\to X \to A_\alpha\to B_\alpha\to0$.
Next we construct $\bE_{\alpha+1}$ in the following steps.

{\bf Step 1.} Construct $A_{\alpha+1}$ and $B_{\alpha+1}$.

Let $i\in V_{\alpha+1}\backslash V_\alpha$.
Then $\oplus_{a\in Q_{1}^{\ast\to i}}B_\alpha(s(a))\in\widetilde{\sf{C}}$,
as each arrow $a\in Q_{1}^{\ast\to i}$ has the source $s(a)$ in $V_\alpha$;
see \ref{rooted}.
Since the cotorsion pair $(\widetilde{\sf{C}},\sf{F})$ is complete,
there exists a short exact sequence
\begin{center}
$0\to \oplus_{a\in Q_{1}^{\ast\to i}}B_\alpha(s(a))\xra{(\epsilon_a)_{a\in Q_{1}^{\ast\to i}}} D_{\alpha+1}^i\to \overline{B}^i\to0$.
\end{center}
in $\A$ with $D_{\alpha+1}^i\in \widetilde{\sf{C}}\cap\sf{F}=\sf{C}\cap\widetilde{\sf{F}}\subseteq\widetilde{\sf{F}}$ and $\overline{B}^i\in\widetilde{\sf{C}}$.
For each $a\in Q_{1}^{\ast\to i}$, there exists the canonical inclusion
\begin{center}
$(\epsilon_a)_{a\in Q_{1}^{\ast\to i}}\circ \iota^{s(a),i}: \,\,
B_\alpha(s(a))\overset{\iota^{s(a),i}}{\hookrightarrow} \oplus_{a\in Q_{1}^{\ast\to i}}B_\alpha(s(a))
\xra{(\epsilon_a)_{a\in Q_{1}^{\ast\to i}}} D_{\alpha+1}^i$
\end{center}
By assumption, one has $i\in Q_0\backslash V_\alpha$.
Hence $B_\alpha(a):B_\alpha(s(a))\to B_\alpha(i)=B(i)$.
So there is a canonical morphism
\begin{equation*}
%\tag{$\ast$}
\left(\begin{array}{cc}
B_\alpha(a) \\ (\epsilon_a)_{a\in Q_{1}^{\ast\to i}}\circ \iota^{s(a),i} \end{array}\right)
:B_\alpha(s(a))\to B(i)\oplus D_{\alpha+1}^i
\end{equation*}
which induces a monomorphism
\begin{equation}
\tag{$\dag$}
\left(\begin{array}{cc}(B_\alpha(a))_{a\in Q_{1}^{\ast\to i}} \\ (\epsilon_a)_{a\in Q_{1}^{\ast\to i}} \end{array}\right)
:\oplus_{a\in Q_{1}^{\ast\to i}}B_\alpha(s(a))\to B(i)\oplus D_{\alpha+1}^i.
\end{equation}

Define $A_{\alpha+1}$ and $B_{\alpha+1}$ as follows:
\begin{equation*}
  A_{\alpha+1}(i)=
  \begin{cases}
    A_{\alpha}(i)                & \text{if $i\notin V_{\alpha+1}\backslash V_\alpha$}\,,\\
    A(i)\oplus D_{\alpha+1}^i    & \text{otherwise}.
  \end{cases}
\end{equation*}
and
\begin{equation*}
  B_{\alpha+1}(i)=
  \begin{cases}
    B_{\alpha}(i)                & \text{if $i\notin V_{\alpha+1}\backslash V_\alpha$}\,,\\
    B(i)\oplus D_{\alpha+1}^i    & \text{otherwise}.
  \end{cases}
\end{equation*}
For an arrow $a: j\to k$ in $Q_1$, the morphisms $A_{\alpha+1}(a)$ and $B_{\alpha+1}(a)$ are defined as:
\begin{rqm}
\item[--] If $k\in V_{\alpha+1}\backslash V_\alpha$, then $j\in V_\alpha$.
          Define
          $$A_{\alpha+1}(a)= \left(\begin{array}{cc}A_\alpha(a) \\ (\epsilon_a)_{a\in Q_{1}^{\ast\to k}}\circ \iota^{j,k}\circ h_\alpha(j)
          \end{array}\right) :A_{\alpha+1}(j)\to A_{\alpha+1}(k)= A(k)\oplus D_{\alpha+1}^k,$$
          and define
          $$B_{\alpha+1}(a)= \left(\begin{array}{cc}B_\alpha(a) \\ (\epsilon_a)_{a\in Q_{1}^{\ast\to k}}\circ \iota^{j,k}
          \end{array}\right) :B_{\alpha+1}(j)\to B_{\alpha+1}(k)= B(k)\oplus D_{\alpha+1}^k.$$

\item[--] If $j\in V_{\alpha+1}\backslash V_\alpha$, then $k\notin V_{\alpha+1}$.
          Hence $A_{\alpha+1}(j)= A(j)\oplus D^{\alpha+1}_j$,
          $A_{\alpha+1}(k)= A_{\alpha}(k)$,
          $B_{\alpha+1}(j)= B(j)\oplus D_{\alpha+1}^j$ and
          $B_{\alpha+1}(k)= B_{\alpha}(k)$.
          Define $A_{\alpha+1}(a)$ as the composition of the projection map $A(j)\oplus D_{\alpha+1}^j\twoheadrightarrow A(j)$
          followed by $A(a)(=A_\alpha(a))$, and define
          $B_{\alpha+1}(a)$ as the composition of the projection map $B(j)\oplus D_{\alpha+1}^j\twoheadrightarrow B(j)$
          followed by $B(a)(=B_\alpha(a))$.

\item[--] For the other cases, define $A_{\alpha+1}(a)=A_{\alpha}(a)$ and $B_{\alpha+1}(a)=B_{\alpha}(a)$.
\end{rqm}
Note that $D_{\alpha+1}^i\in \widetilde{\sf{C}}\cap\sf{F}=\sf{C}\cap\widetilde{\sf{F}}\subseteq\widetilde{\sf{F}}$
for each $i\in V_{\alpha+1}\backslash V_\alpha$.
It follows that $A(i)\oplus D_{\alpha+1}^i\in\widetilde{\sf{F}}$.
Thus $A_{\alpha+1}(i)\in \widetilde{\sf{F}}$ for all $i\in Q_0$.
On the other hand, for each $i\in V_{\alpha+1}\backslash V_\alpha$,
$B(i)\oplus D_{\alpha+1}^i\in \widetilde{\sf{C}}$.
Hence $B_{\alpha+1}(i)\in \widetilde{\sf{C}}$ for all $i\in Q_0$.
Thus, to show $A_{\alpha+1}$ and $B_{\alpha+1}$ constructed above
satisfy the desired conditions (b) and (c),
it remains to prove that $\C_i(B_{\alpha+1})\in \widetilde{\sf{C}}$
for all $i\in V_{\alpha+1}\backslash V_\alpha$,
that is, the cokernel of the morphism
$$\varphi^{B_{\alpha+1}}_i=
\left(\begin{array}{cc}(B_\alpha(a))_{a\in Q_{1}^{\ast\to i}}\\ (\epsilon_a)_{a\in Q_{1}^{\ast\to i}}\end{array}\right)$$
given in $(\dag)$ is in $\widetilde{\sf{C}}$.
Consider the following diagram with exact rows and columns
$$\xymatrix@C=18pt@R=15pt{
 & &0 \ar[d]    \\
&& B(i)\ar[d]^{\iota}  \\
0 \ar[r] & \oplus_{a\in Q_{1}^{\ast\to i}}B_\alpha(s(a)) \ar[r]^{\ \ \ \ \varphi^{B_{\alpha+1}}_i}\ar@{=}[d] & B(i)\oplus D_{\alpha+1}^i \ar[r]\ar[d]^{\pi} & \C_i(B_{\alpha+1})     \ar[r]\ar@.[d] & 0 \\
0 \ar[r] & \oplus_{a\in Q_{1}^{\ast\to i}}B_\alpha(s(a)) \ar[r]^{} &  D_{\alpha+1}^i  \ar[r]\ar[d] & \overline{B}^i     \ar[r]& 0. \\
& &0  }
$$
The left square is commutative clearly, so there exists a morphism $\C_i(B_{\alpha+1})\to \overline{B}^i$ such that the right square is commutative. By the Snake Lemma, one gets an exact sequence
$0\to B(i)\to \C_i(B_{\alpha+1})\to \overline{B}^i\to 0$ in $\A$.
Since both $B(i)$ and $\overline{B}^i$ are in $\widetilde{\sf{C}}$,
we conclude that $\C_i(B_{\alpha+1})\in \widetilde{\sf{C}}$ as well.

{\bf Step 2.} Construct $\bE_{\alpha+1}$.

Set $X_{\alpha+1}=X$, and define $k_{\alpha+1}$ and $h_{\alpha+1}$ as follows:
\begin{equation*}
  k_{\alpha+1}(i)=
  \begin{cases}
    k_{\alpha}(i)               & \text{if $i\notin V_{\alpha+1}\backslash V_\alpha$}\,,\\
   \\
    \left(
                       \begin{array}{cc}
                         k_{\alpha}(i) \\
                         0 \\
                       \end{array}
                     \right)    & \text{otherwise}.
  \end{cases}
\end{equation*}
and
\begin{equation*}
 \quad \, h_{\alpha+1}(i)=
  \begin{cases}
    h_{\alpha}(i)               & \text{if $i\notin V_{\alpha+1}\backslash V_\alpha$}\,,\\
   \\
   \left( \begin{array}{cc}
                                     h_{\alpha}(i) & 0 \\
                                                      0 & 1 \\
                                                    \end{array}
                                                  \right) & \text{otherwise}.
  \end{cases}
\end{equation*}
According to the constructions of $A_{\alpha+1}$ and $B_{\alpha+1}$
and combining the constructions of $k_{\alpha+1}$ and $h_{\alpha+1}$,
one gets the desired exact sequence
$$\bE_{\alpha+1}:\,\,\,
0\to X_{\alpha+1} \xra{k_{\alpha+1}} A_{\alpha+1}\xra{h_{\alpha+1}} B_{\alpha+1}\to0.$$

{\bf Step 3.} Construct monomorphisms $f_{{\alpha+1},\alpha}$ and $g_{{\alpha+1},\alpha}$.

Define $f_{{\alpha+1},\alpha}$ and $g_{{\alpha+1},\alpha}$ as follows:
\begin{equation*}
  f_{{\alpha+1},\alpha}(i)=
  \begin{cases}
    \textrm{id}_{A_\alpha(i)}              & \text{if $i\notin V_{\alpha+1}\backslash V_\alpha$}\,,\\
   \\
    \left(
                       \begin{array}{cc}
                         1 \\
                         0 \\
                       \end{array}
                     \right)    & \text{otherwise}.
  \end{cases}
\end{equation*}
and
\begin{equation*}
g_{{\alpha+1},\alpha}(i)=
  \begin{cases}
    \textrm{id}_{B_\alpha(i)}              & \text{if $i\notin V_{\alpha+1}\backslash V_\alpha$}\,,\\
   \\
    \left(
                       \begin{array}{cc}
                         1 \\
                         0 \\
                       \end{array}
                     \right)    & \text{otherwise}.
  \end{cases}
\end{equation*}
It is clear that $f_{{\alpha+1},\alpha}$ and $g_{{\alpha+1},\alpha}$ are monomorphisms, and the next diagram is commutative
$$\xymatrix{
\bE_\alpha:\quad 0\ar[r] &X\ar[r]^{\quad k_\alpha\ \ \ }\ar@{=}[d] & A_\alpha \ar[d]_{f_{{\alpha+1},\alpha}} \ar[r]^{h_\alpha}
&  B_\alpha \ar[d]_{{g_{{\alpha+1},\alpha}}} \ar[r]&0  \\
\bE_{\alpha+1}:\quad 0\ar[r] &X\ar[r]^{\quad k_{\alpha+1}\ \ \ \ \ }&  A_{\alpha+1}  \ar[r]^{h_{\alpha+1}} & B_{\alpha+1}\ar[r]&0.}$$

Suppose now that $\beta\leqslant\lambda$ is a limit ordinal and $\bE_\alpha$ is constructed for all $\alpha<\beta$.
Next we construct $\bE_\beta$.
In this case one has $V_\beta=\cup_{\alpha<\beta}V_{\alpha}$.
\begin{rqm}
\item[--] If $i\in V_\beta$, then $i\in V_\alpha$ for some ordinal $\alpha<\beta$,
          and so for all ordinal $\alpha<\alpha'\leqslant\beta$
          one has $\bE_{\alpha'}(i)=\bE_{\alpha}(i)$ as $i\notin V_{\alpha'}\backslash V_{\alpha}$.
\item[--] If $i\notin V_\beta$, then by hypothesis, one has $\bE_{\alpha}(i)=\bE(i)$ for all $\alpha<\beta$.
\end{rqm}
We let $\bE_\beta=\colim_{\alpha<\beta}\bE_\alpha$. Then one has
\begin{equation*}
  \bE_{\beta}(i)= \colim_{\alpha<\beta}\bE_\alpha(i)=
  \begin{cases}
    \bE_{\alpha}(i)\ \text{for some~} \alpha<\beta  & \text{if $i\in V_{\beta}$}\,,\\
    \bE(i) & \text{if $i\notin V_{\beta}$}\,.
  \end{cases}
\end{equation*}
Thus $\bE_{\beta}$ is a short exact sequence in $\Rep{Q}{\A}$ satisfying the desired conditions.

Finally, let $A'=A_\lambda$ and $B'=B_\lambda$.
Then one gets the desired exact sequence $\bE'=\bE_\lambda:\,\,\,0\to X\to A'\to B'\to 0$.
\end{prf*}

Now Theorem B in the introduction holds by Lemmas \ref{model structor} and \ref{trivial objects}. It is known that $(\GF{R},\PGF{R}^\perp,\Cot{R})$ is a hereditary Hovey triple on $\Mod{R}$ (see \cite[pp. 27]{SS20}), so Theorem C in the introduction follows from Theorems A and B, Lemma \ref{cotorsion} and Corollary \ref{pgfpair}.

\begin{rmk}
Dually, one can prove that if $Q$ is a right rooted quiver and
$\sf{A}$ is an abelian category with enough projectives
which satisfies the axiom $\sf{AB4}^\ast$ then any hereditary Hovey triple $\sf{(C,W,F)}$ on $\sf{A}$
induces an hereditary Hovey triple
$$(\Rep{Q}{\sf{C}},\Rep{Q}{\W},\Psi(\sf{F}))$$
on $\Rep{Q}{\sf{A}}$; we leave the proof to the reader.
\end{rmk}

\appendix
\section*{Appendix}
\stepcounter{section}
\noindent
In this section we reprove Lemmas \ref{cogenerator} and \ref{trivial objects} over a simple left rooted quiver to comprehend the ideas of the proofs. Let $Q$ be the quiver
$$
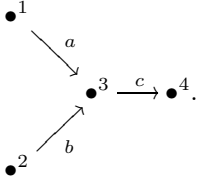
\xymatrix@C=15pt@R=15pt{
  \bullet^{1} \ar[dr]^{a} \\
  & \bullet^{3} \ar[r]^{c}  &  \bullet^{4}. \\
  \bullet^{2}  \ar[ur]_{b} }$$
Then $Q$ is a left rooted quiver with $V_0=\emptyset$, $V_1=\{1,2\}$, $V_2=\{1,2,3\}$, $V_3=Q_0$.

\begin{bfhpg}[\bf Proof of Lemma \ref{cogenerator}]\label{pf1}
Let $X$ be in $\Phi(\sf X)$. Next we construct an exact sequence $0 \to X \to W \to Y \to 0$ in $\Rep{Q}{R}$ with $W\in\Phi(\sf W)$ and $Y\in\Phi(\sf X)$.

We let $\bE_0$ be the short exact sequence of zero representations.

Since $\sf W$ is a cogenerator for $\sf X$,
we conclude that there exist exact sequences $0 \to X(1) \to W^1 \to Y^1 \to 0$ and $0 \to X(2) \to W^2 \to Y^2 \to 0$ with $W^1,W^2\in\sf W$ and $Y^1,Y^2\in\sf X$. So we let $\bE_1: 0 \to X_1 \to W_1 \to Y_1 \to 0$ be the exact sequence in $\Rep{Q}{R}$ as follows:
\begin{center}
\setlength{\unitlength}{0.1in}
\begin{picture}(45,5.5)
\put(0,0){$0 \to
\begin{pmatrix}
       X(1) &          \\
            & \ 0 \to 0  \\
       X(2) &
\end{pmatrix}
\to
\begin{pmatrix}
       W^1 &          \\
           & \ 0 \to 0  \\
       W^2 &
\end{pmatrix}
\to
\begin{pmatrix}
       Y^1 &          \\
            & \ 0 \to 0  \\
       Y^2 &
\end{pmatrix}
\to 0 \:.$}
\put(7.2,1){$\searrow$}\put(7.2,-1){$\nearrow$}%\put(9,-1.5){$f$}\put(9,1.5){$f$}
\put(19.2,1){$\searrow$}\put(19.2,-1){$\nearrow$}
\put(31,1){$\searrow$}\put(31,-1){$\nearrow$}
%\put(0,0){\line(1,0){45}}这是一条测试用的水平参考直线
%\put(0,-2.5){\line(0,1){5.5}}这是一条测试用的垂直参考直线
\end{picture}
\end{center}
\vspace{1cm}
Clearly, for all $i\in V_1$, $\varphi_{i}^{W_1}$ and $\varphi_{i}^{Y_1}$ are injective with $\C_i(W_1)\in\sf W$ and $\C_i(Y_1)\in\sf X$.

Since $(\sf W, {\sf W}^{\perp})$ is a complete cotorsion pair, there is an exact sequence $0\to W^1\oplus W^2 \to U^3 \to C' \to 0$ with $U^3\in{\sf W}^{\perp}$ and $C'\in\sf W$. It is clear that $W^1\oplus W^2$ is in $\sf W$, so one has $U^3\in\sf W\cap{\sf W}^{\perp}$. Consider the following push-out diagram:
$$\xymatrix@C=10pt@R=20pt{&&0\ar[d] &0\ar[d]\\
  0 \ar[r] & X(1)\oplus X(2) \ar[r]^{}\ar@{=}[d] & W^1\oplus W^2 \ar[r]\ar[d]^{\iota} & Y^1\oplus Y^2 \ar[r]\ar[d]^{\tau} & 0\\
  0 \ar[r] & X(1)\oplus X(2) \ar[r] & U^3 \ar[r]^{}\ar[d] & T^3 \ar[r]\ar[d] & 0.\\
  && C' \ar@{=}[r]\ar[d] & C' \ar[d]\\
  && 0 & 0 }$$
Since $\sf W$ is a cogenerator for $\sf X$ and $\C_3(X)\in\sf X$, one has an exact sequence
$$0 \to \C_3(X) \to S^3 \to Z^3 \to 0$$
with $S^3 \in \sf W$ and $Z^3\in\sf X$. By the assumption one has $\Ext{1}{\C_3(X)}{U^3}=0$. Thus by Lemma \ref{diagram}, one gets a commutative diagram with exact rows and columns
$$\xymatrix@C=15pt@R=15pt{
  & 0 \ar[d] & 0 \ar[d] & 0 \ar[d] \\
  0 \ar[r] & X(1)\oplus X(2) \ar[d]\ar[r]^{\ \ \ \ \ \ \varphi_3^{X}} & X(3) \ar[d]\ar[r] & \C_3(X) \ar[d]\ar[r] & 0 \\
  0 \ar[r] & U^3 \ar[d]\ar[r]^{\iota'} & W^3 \ar[d]\ar[r] & S^3 \ar[d]\ar[r] & 0 \\
  0 \ar[r] & T^3 \ar[d]\ar[r]^{\tau'} & Y^3 \ar[d]\ar[r] &Z^3 \ar[d]\ar[r] & 0, \\
  & 0 & 0 & 0  }$$
where $W^3=U^3\oplus S^3\in\sf W$ and $Y^3\in\sf X$, as $\sf W$ and $\sf X$ are closed under extensions.
We let $\bE_2: 0 \to X_2 \to W_2 \to Y_2 \to 0$ be the exact sequence in $\Rep{Q}{R}$ as follows:
\begin{center}
\setlength{\unitlength}{0.1in}
\begin{picture}(45,5.5)
\put(0,0){$0 \to
\begin{pmatrix}
       X(1) &          \\
            & \ X(3) \to 0  \\
       X(2) &
\end{pmatrix}
\to
\begin{pmatrix}
       W^1 &          \\
           & \ W^3 \to 0  \\
       W^2 &
\end{pmatrix}
\to
\begin{pmatrix}
       Y^1 &          \\
            & \ Y^3 \to 0  \\
       Y^2 &
\end{pmatrix}
\to 0 \:.$}
\put(7.2,1){$\searrow$}\put(7.2,-1){$\nearrow$}%\put(9,-1.5){$f$}\put(9,1.5){$f$}
\put(21.6,1){$\searrow$}\put(21.6,-1){$\nearrow$}\put(22.4,-1.5){$\iota'\iota\epsilon_2$}\put(22.2,1.8){$\iota'\iota\epsilon_1$}
\put(34.5,1){$\searrow$}\put(34.5,-1){$\nearrow$}\put(35.4,-1.5){$\tau'\tau\epsilon_2$}\put(35.4,1.8){$\tau'\tau\epsilon_1$}
%\put(0,0){\line(1,0){45}}这是一条测试用的水平参考直线
%\put(0,-2.5){\line(0,1){5.5}}这是一条测试用的垂直参考直线
\end{picture}
\end{center}
\vspace{1cm}
where $\epsilon_i$ is the natural embedding. The homomorphisms $\varphi_{i}^{W_2}$ and $\varphi_{i}^{Y_2}$ are injective with $\C_i(W_2)\in\sf W$ and $\C_i(Y_2)\in\sf X$ for $i=1,2$. It is easy to see that $\varphi_{3}^{W_2}=\iota'\iota$ and $\varphi_{3}^{Y_2}=\tau'\tau$ are injective. Next we show that $\C_3(W_2)\in\sf W$ and $\C_3(Y_2)\in\sf X$; we only prove the case for $W_2$. Using the Snake Lemma, one gets the next commutative diagram with exact rows and columns:
$$\xymatrix@C=20pt@R=20pt{&&0\ar[d] &0\ar[d]\\
  0 \ar[r] & W^1\oplus W^2 \ar[r]\ar@{=}[d] & U^3 \ar[r]\ar[d] & C' \ar[r]\ar[d] & 0\\
  0 \ar[r] & W^1\oplus W^2 \ar[r]^{\ \ \ \ \ \ \varphi_3^{W_{\alpha+1}}} & W^3 \ar[r]\ar[d] & \C_3(W_{2}) \ar[r]\ar[d] & 0.\\
  && S^3 \ar@{=}[r]\ar[d] & S^3 \ar[d]\\
  && 0 & 0 }$$
Since $C'$ and $S^3$ are in $\sf W$, so is $\C_3(W_{2})$ as $\sf W$ is closed under extensions.

Finally, $\bE_3$ is constructed similarly; we still give its construction here. Fix an exact sequence $0\to W^3 \to U^4 \to C'' \to 0$ with $U^4\in\sf W\cap{\sf W}^{\perp}$ and $C''\in\sf W$. Consider the following push-out diagram:
$$\xymatrix@C=15pt@R=20pt{&&0\ar[d] &0\ar[d]\\
  0 \ar[r] & X(3) \ar[r]\ar@{=}[d] & W^3 \ar[r]\ar[d]^{\mu} & Y^3 \ar[r]\ar[d]^{\nu} & 0\\
  0 \ar[r] & X(3) \ar[r] & U^4 \ar[r]^{}\ar[d] & T^4 \ar[r]\ar[d] & 0.\\
  && C'' \ar@{=}[r]\ar[d] & C'' \ar[d]\\
  && 0 & 0 }$$
Since $\sf W$ is a cogenerator for $\sf X$ and $\C_4(X)\in\sf X$, one has an exact sequence
$$0 \to \C_4(X) \to S^4 \to Z^4 \to 0$$
with $S^4 \in \sf W$ and $Z^4\in\sf X$. By the assumption one has $\Ext{1}{\C_4(X)}{U^4}=0$. Thus by Lemma \ref{diagram}, one gets a commutative diagram with exact rows and columns
$$\xymatrix@C=15pt@R=15pt{
  & 0 \ar[d] & 0 \ar[d] & 0 \ar[d] \\
  0 \ar[r] & X(3) \ar[d]\ar[r]^{\ \varphi_4^{X}} & X(4) \ar[d]\ar[r] & \C_4(X) \ar[d]\ar[r] & 0 \\
  0 \ar[r] & U^4 \ar[d]\ar[r]^{\mu'} & W^4 \ar[d]\ar[r] & S^4 \ar[d]\ar[r] & 0 \\
  0 \ar[r] & T^4 \ar[d]\ar[r]^{\nu'} & Y^4 \ar[d]\ar[r] &Z^4 \ar[d]\ar[r] & 0, \\
  & 0 & 0 & 0  }$$
where $W^4=U^4\oplus S^4\in\sf W$ and $Y^4\in\sf X$, as $\sf W$ and $\sf X$ are closed under extensions.
We let $\bE_3: 0 \to X_3 \to W_3 \to Y_3 \to 0$ be the exact sequence in $\Rep{Q}{R}$ as follows:
\setlength{\unitlength}{0.1in}
\begin{picture}(45,5.5)
\put(0,0){$0 \to
\begin{pmatrix}
       X(1) &          \\
            & \ X(3) \to X(4)  \\
       X(2) &
\end{pmatrix}
\to
\begin{pmatrix}
       W^1 &          \\
           & \ W^3 \to W^4  \\
       W^2 &
\end{pmatrix}
\to
\begin{pmatrix}
       Y^1 &          \\
            & \ Y^3 \to Y^4  \\
       Y^2 &
\end{pmatrix}
\to 0 \:.$}
\put(7.2,1){$\searrow$}\put(7.2,-1){$\nearrow$}
\put(23.6,1){$\searrow$}\put(23.6,-1){$\nearrow$}
\put(24.4,-1.5){$\iota'\iota\epsilon_2$}\put(24.2,1.8){$\iota'\iota\epsilon_1$}\put(27.8,1){$\mu'\mu$}
\put(38.5,1){$\searrow$}\put(38.5,-1){$\nearrow$}
\put(39.5,-1.5){$\tau'\tau\epsilon_2$}\put(39,1.8){$\tau'\tau\epsilon_1$}\put(42,0.9){$\upsilon'\upsilon$}
%\put(0,0){\line(1,0){45}}这是一条测试用的水平参考直线
%\put(0,-2.5){\line(0,1){5.5}}这是一条测试用的垂直参考直线
\end{picture}
\vspace{1cm}

\noindent
As shown above, the homomorphisms $\varphi_{i}^{W_3}$ and $\varphi_{i}^{Y_3}$ are injective with $\C_i(W_3)\in\sf W$ and $\C_i(Y_3)\in\sf X$ for $i=1,2,3$. It is clear that $\varphi_{4}^{W_3}=\mu'\mu$ and $\varphi_{4}^{Y_3}=\nu'\nu$ are injective. The same argument as above yields that $\C_4(W_3)\in\sf W$ and $\C_4(Y_3)\in\sf X$. Thus one has $W_3\in\Phi(\sf W)$ and $Y_3\in\Phi(\sf X)$. Note that $X_3$ is actually the representation $X$. Then $\Phi(\sf W)$ is a cogenerator for $\Phi(\sf X)$. \qed
\end{bfhpg}

Next, we give the proof of Lemma \ref{trivial objects} over the quiver $Q$ given at the begin of this section.

\begin{bfhpg}[\bf Proof of Lemma \ref{trivial objects}]\label{pf2}
We only explain the inclusion ``$\Rep{Q}{\W}\subseteq\mathsf{T}$".
Let $X\in \Rep{Q}{\W}$.
Then for each $i\in \{1,2,3,4\}$,
there exists a short exact sequence
$$0\to X(i) \overset{k(i)}\longrightarrow A(i) \overset{h(i)}\longrightarrow B(i) \to 0$$
in $\A$ with $A(i)\in \widetilde{\sf{F}}$ and $B(i)\in \widetilde{\sf{C}}$ (see Proposition \ref{model structor111}).
Note that $\widetilde{\sf{F}}\subseteq\sf{F}$.
It follows that for each arrow $e: i\to j \in \{a,b,c\}$,
there exist morphisms $A(e)$ and $B(e)$ such that the following diagram is commutative
$$\xymatrix{
0\ar[r] &X(i)\ar[r]^{k(i)}\ar[d]_{X(e)}& A(i) \ar@.[d]_{A(e)} \ar[r]^{h(i)}&  B(i) \ar@.[d]_{B(e)}\ar[r]&0  \\
0\ar[r] &X(j)\ar[r]^{k(j)}&  A(j)  \ar[r]^{h(j)} & B(j)\ar[r]&0}$$
Hence there is a short exact sequence
$$\bE:\,\,\,0\to X\to A\to B\to 0$$
in $\Rep{Q}{\A}$ with the following form:

$$\xymatrix@R=0.01cm{
  &&& 0 \ar[ddd] \\
  &0\ar[ddd]&&&& 0 \ar[ddd] &  0 \ar[ddd]\\
  &&&& 0 \ar[ddd]\\
  &&& X(1)\ar[drr]^{\quad \quad X(a)}\ar[ddd]_{k(1)} \\
  &~X:\ar[ddd]&&&& X(3) \ar[r]^{X(c)} \ar[ddd]_{k(3)} &  X(4)\ar[ddd]_{k(4)} \\
  &&&& X(2)  \ar[ur]_{X(b)} \ar[ddd]_{k(2)}\\
  &&& A(1) \ar@.[drr]^{\quad \quad A(a)} \ar[ddd]_{h(1)}\\
  &~A:\ar[ddd]&&&& A(3) \ar@.[r]^{A(c)} \ar[ddd]_{h(3)} &  A(4)\ar[ddd]_{h(4)} \\
  &&&& A(2)  \ar@.[ur]_{A(b)}\ar[ddd]_{h(2)}\\
  &&& B(1) \ar@.[drr]^{\quad \quad B(a)} \ar[ddd]\\
  &~B:\ar[ddd]&&&& B(3) \ar@.[r]^{B(c)} \ar[ddd] &  B(4)\ar[ddd]\\
  &&&& B(2)  \ar@.[ur]_{B(b)}\ar[ddd]\\
  &&&0  \\
  &0&&&& 0  &  0 \\
  &&&&0
  }$$
It is clear that $A\in \Rep{Q}{\widetilde{\sf{F}}}$ and $B\in \Rep{Q}{\widetilde{\sf{C}}}$, and
$\varphi^B_1$ and $\varphi^B_2$ are monomorphisms
with $\C_1(B)=B(1)$ and $\C_2(B)=B(2)$ in $\widetilde{\sf{C}}$,
as the vertices 1 and 2 are in $V_1$.
So we let $\bE_1=\bE$.
But $\varphi^B_3$ and $\varphi^B_4$ may fail to be monomorphisms,
and $\C_3(B)$ and $\C_4(B)$ may fail to be in $\widetilde{\sf{C}}$.
Next we repair $\bE_1$ at the vertex 3.

Since $B(1)\oplus B(2)\in\widetilde{\sf{C}}$ and the cotorsion pair $(\widetilde{\sf{C}}, \sf{F})$ is complete,
there exists a short exact sequence
\begin{center}
$0\to B(1)\oplus B(2)\xra{(\epsilon_2^1,\epsilon_2^2)} D_{2}^3\to {\overline{B}^3}\to0$.
\end{center}
in $\A$ with $D_{2}^3\in \widetilde{\sf{C}}\cap\sf{F}=\sf{C}\cap\widetilde{\sf{F}}\subseteq\widetilde{\sf{F}}$ and ${\overline{B}}^3\in\widetilde{\sf{C}}$.
Now we let $A_2$ and $B_2$ be representations as follows:
$$\xymatrix@R=0.5cm{
  &&A(1) \ar[dr]^{A_2(a)} \\
  &A_2:&& A(3)\oplus D_2^3 \ar[r]^{\quad A_2(c)}  &  A(4) \\
  &&A(2)  \ar[ur]_{A_2(b)} }$$
and
$$\xymatrix@R=0.5cm{
  &&B(1) \ar[dr]^{B_2(a)} \\
   &B_2:&& B(3)\oplus D_2^3 \ar[r]^{\quad B_2(c)}  &  B(4) \\
  &&B(2)  \ar[ur]_{B_2(b)} }$$
where $A_2(a)= \left(\begin{array}{cc}A(a) \\
                                    (\epsilon_2^1,\epsilon_2^2)\circ \iota\circ h(1)
                       \end{array}\right)$,
$A_2(b)= \left(\begin{array}{cc}A(b) \\
                                    (\epsilon_2^1,\epsilon_2^2)\circ \iota\circ h(2)
                       \end{array}\right)$,
$A_2(c)=A(c)\circ\pi$,
$B_2(a)= \left(\begin{array}{cc}B(a) \\
                                    (\epsilon_2^1,\epsilon_2^2)\circ \iota
                       \end{array}\right)$,
$B_2(b)= \left(\begin{array}{cc}B(b) \\
                                    (\epsilon_2^1,\epsilon_2^2)\circ \iota
                       \end{array}\right)$
and $B_2(c)=B(c)\circ\pi$
(here $\iota$ and $\pi$ are the canonical inclusion and projection morphisms, respectively).

Note that $D_{2}^3\in\widetilde{\sf{F}}$.
It follows that $A(3)\oplus D_{2}^3\in\widetilde{\sf{F}}$.
Thus $A_{2}(i)\in \widetilde{\sf{F}}$ for all $i\in \{1,2,3,4\}$.
This implies that $A_2\in \Rep{Q}{\widetilde{\sf{F}}}$.
On the other hand, since $B(3)\oplus D_{2}^3\in \widetilde{\sf{C}}$,
we conclude that $B_{2}(i)\in \widetilde{\sf{C}}$ for all $i\in \{1,2,3,4\}$.
It is easy to see that the morphism
$$\varphi^{B_2}_3=
\left(\begin{array}{cc}
(B(a), B(b))\\
(\epsilon_2^1,\epsilon_2^2)
      \end{array}\right)$$
is now a monomorphism.
Next we prove that $\C_3(B_{2})= \Coker\,\varphi^{B_2}_3\in \widetilde{\sf{C}}$.
To this end, consider the following diagram with exact rows and columns
$$\xymatrix@C=15pt@R=15pt{
 & &0 \ar[d]    \\
&& B(3)\ar[d]^{\iota}  \\
0 \ar[r] & B(1)\oplus B(2) \ar[r]^{\ \ \varphi^{B_2}_3}\ar@{=}[d] & B(3)\oplus D_{2}^3 \ar[r]\ar[d]^{\pi} & \C_3(B_{2})     \ar[r]\ar@.[d] & 0 \\
0 \ar[r] & B(1)\oplus B(2) \ar[r]^{} &  D_{2}^3  \ar[r]\ar[d] & \overline{B}^3     \ar[r]& 0. \\
& &0  }
$$
The left square is commutative, so there exists a morphism $\C_3(B_{2})\to \overline{B}^3$
such that the right square is commutative. By the Snake Lemma, one gets an exact sequence
$0\to B(3)\to \C_3(B_{2})\to \overline{B}^3\to 0$
in $\A$.
Since both $B(3)$ and $\overline{B}^3$ are in $\widetilde{\sf{C}}$,
we conclude that $\C_3(B_{2})$ is in $\widetilde{\sf{C}}$ as well.

Define $k_{2}$ and $h_{2}$ as follows:
\begin{equation*}
  k_{2}(i)=
  \begin{cases}
    k_1(i)=k(i)               & \text{if $i\in\{1,2,4\}$}\,,\\
   \\
    \left(
                       \begin{array}{cc}
                          k_1(3)=k(3) \\
                         0 \\
                       \end{array}
                     \right)    & \text{if $i=3$}.
  \end{cases}
\end{equation*}
and
\begin{equation*}
 \quad \, h_{2}(i)=
  \begin{cases}
    h_1(i)=h(i)               & \text{if $i\in\{1,2,4\}$}\,,\\
   \\
   \left( \begin{array}{cc}
                                     h_1(3)=h(3) & 0 \\
                                                      0 & 1 \\
                                                    \end{array}
                                                  \right) & \text{if $i=3$}.
  \end{cases}
\end{equation*}
According to the constructions of $A_{2}$ and $B_{2}$
and combining the constructions of $k_{2}$ and $h_{2}$,
one gets an exact sequence
$$\bE_{2}:\,\,\,0\to X \overset{k_{2}}\longrightarrow A_{2}\overset{h_{2}}\longrightarrow B_{2}\to0$$
in $\Rep{Q}{\A}$ with the following form

$$\xymatrix@R=0.01cm{
  &&& 0 \ar[ddd] \\
  &0\ar[ddd]&&&& 0 \ar[ddd] &  0 \ar[ddd]\\
  &&&& 0 \ar[ddd]\\
  &&& X(1)\ar[drr]^{\quad \quad X(a)}\ar[ddd]_{k(1)} \\
  &~X:\ar[ddd]_{k_{2}}&&&& X(3) \ar[r]^{X(c)} \ar[ddd]_{k_2(3)} &  X(4)\ar[ddd]_{k_2(4)} \\
  &&&& X(2)  \ar[ur]_{X(b)} \ar[ddd]_{ k(2)}\\
  &&& A(1) \ar[drr]^{\quad \quad A_2(a)} \ar[ddd]_{h(1)}\\
  &~A_2:\ar[ddd]_{h_2}&&&& A(3)\oplus D_2^3 \ar[r]^{\quad A_2(c)} \ar[ddd]_{h_2(3)} &  A(4)\ar[ddd]_{h_2(4)} \\
  &&&& A(2)  \ar[ur]_{A_2(b)\quad}\ar[ddd]_{h(2)}\\
  &&& B(1) \ar[drr]^{\quad \quad B_2(a)} \ar[ddd]\\
  &~B_2:\ar[ddd]&&&& B(3)\oplus D_2^3 \ar[r]^{\quad B_2(c)} \ar[ddd] &  B(4)\ar[ddd]\\
  &&&& B(2)  \ar[ur]_{B_2(b)}\ar[ddd]\\
  &&&0  \\
  &0&&&& 0  &  0 \\
  &&&&0}$$
Now, all $\varphi^{B_2}_1$, $\varphi^{B_2}_2$ and $\varphi^{B_2}_3$ are monomorphisms
with all $\C_1(B_2)=B(1)$, $\C_2(B_2)=B(2)$ and $\C_3(B_2)$ in $\widetilde{\sf{C}}$.
But $\varphi^{B_2}_4$ may fail to be a monomorphism,
and $\C_4(B_2)$ may fail to be in $\widetilde{\sf{C}}$.
We proceed the same argument for repairing $\bE_{2}$ at the vertex 4.

Since $B(3)\oplus D_2^3\in\widetilde{\sf{C}}$ and the cotorsion pair $(\widetilde{\sf{C}}, \sf{F})$ is complete,
there exists a short exact sequence
\begin{center}
$0\to B(3)\oplus D_2^3\xra{(\epsilon_3^1,\epsilon_3^2)} D_{3}^4\to {\overline{B}^4}\to0$.
\end{center}
in $\A$ with $D_{3}^4\in \widetilde{\sf{C}}\cap\sf{F}=\sf{C}\cap\widetilde{\sf{F}}\subseteq\widetilde{\sf{F}}$ and ${\overline{B}}^4\in\widetilde{\sf{C}}$.
Let $A_3$ and $B_3$ be the representations as follows:
$$\xymatrix@R=0.5cm{
  && A(1) \ar[dr]^{A_3(a)} \\
  &A_3:&& A(3)\oplus D_2^3 \ar[r]^{A_3(c)}  &  A(4)\oplus D_{3}^4 \\
  &&A(2)  \ar[ur]_{A_3(b)} }$$
and
$$\xymatrix@R=0.5cm{
  &&B(1) \ar[dr]^{B_3(a)} \\
   &B_3:&& B(3)\oplus D_2^3 \ar[r]^{B_3(c)}  &  B(4)\oplus D_{3}^4 \\
  &&B(2)  \ar[ur]_{B_3(b)} }$$
where
$A_3(a)= A_2(a)$,
$A_3(b)= A_2(b)$,
$A_3(c)= \left(
           \begin{array}{cc}
              A(b)                  &      0       \\
             \epsilon_3^1\circ h(3) & \epsilon_3^2 \\
           \end{array}
         \right)$,
$B_3(a)= B_2(a)$,
$B_3(b)= B_2(b)$ and
$B_3(c)=\left(
           \begin{array}{cc}
              B(b)                  &      0       \\
             \epsilon_3^1           & \epsilon_3^2 \\
           \end{array}
         \right)$.

Note that $D_{3}^4\in\widetilde{\sf{F}}$.
It follows that $A(4)\oplus D_{3}^4\in\widetilde{\sf{F}}$.
Thus $A_{3}(i)\in \widetilde{\sf{F}}$ for all $i\in \{1,2,3,4\}$.
This implies that $A_3\in \Rep{Q}{\widetilde{\sf{F}}}$.
On the other hand, since $B(4)\oplus D_{3}^4\in \widetilde{\sf{C}}$,
we conclude that $B_{3}(i)\in \widetilde{\sf{C}}$ for all $i\in \{1,2,3,4\}$.
It is easy to see that $\varphi^{B_3}_4=B_3(c)$
is now a monomorphism.
Next we prove that $\C_4(B_{3})= \Coker\,\varphi^{B_3}_4\in \widetilde{\sf{C}}$.
To this end, consider the following diagram with exact rows and columns
$$\xymatrix@C=15pt@R=15pt{
 & &0 \ar[d]    \\
&& B(4)\ar[d]^{\iota}  \\
0 \ar[r] & B(3)\oplus D_2^3 \ar[r]^{\varphi^{B_3}_4}\ar@{=}[d] & B(4)\oplus D_{3}^4 \ar[r]\ar[d]^{\pi} & \C_4(B_{3})     \ar[r]\ar@.[d] & 0 \\
0 \ar[r] & B(3)\oplus D_2^3 \ar[r]^{} &  D_{3}^4  \ar[r]\ar[d] & \overline{B}^4     \ar[r]& 0. \\
& &0  }
$$
The left square is commutative,
so there exists a morphism $\C_4(B_{3})\to \overline{B}^4$
such that the right square is commutative.
By the Snake Lemma, one gets a short exact sequence
$0\to B(4)\to \C_4(B_{3})\to \overline{B}^4\to 0$
in $\A$.
Since both $B(4)$ and $\overline{B}^4$ are in $\widetilde{\sf{C}}$,
we conclude that $\C_4(B_{3})\in \widetilde{\sf{C}}$ as well.

Define $k_{3}$ and $h_{3}$ as follows:
\begin{equation*}
  k_{3}(i)=
  \begin{cases}
    k_2(i)               & \text{if $i\in\{1,2,3\}$}\,,\\
   \\
    \left(
                       \begin{array}{cc}
                         k_2(4) \\
                         0 \\
                       \end{array}
                     \right)    & \text{if $i=4$}.
  \end{cases}
\end{equation*}
and
\begin{equation*}
 \quad \, h_{3}(i)=
  \begin{cases}
    h_2(i)               & \text{if $i\in\{1,2,3\}$}\,,\\
   \\
   \left( \begin{array}{cc}
                                     h_2(4) & 0 \\
                                                      0 & 1 \\
                                                    \end{array}
                                                  \right) & \text{if $i=4$}.
  \end{cases}
\end{equation*}
According to the constructions of $A_{3}$ and $B_{3}$
and combining the constructions of $k_{3}$ and $h_{3}$,
one gets an exact sequence
$$\bE_{3}:\,\,\,0\to X \overset{k_{3}}\longrightarrow A_{3}\overset{h_{3}}\longrightarrow B_{3}\to0$$
in $\Rep{Q}{\A}$ with the following form

$$\xymatrix@R=0.01cm{
  &&& 0 \ar[ddd] \\
  &0\ar[ddd]&&&& 0 \ar[ddd] &  0 \ar[ddd]\\
  &&&& 0 \ar[ddd]\\
  &&& X(1)\ar[drr]^{\quad \quad X(a)}\ar[ddd]_{k(1)} \\
  &~X:\ar[ddd]_{k_{3}}&&&& X(3) \ar[r]^{X(c)} \ar[ddd]_{k_3(3)} &  X(4)\ar[ddd]_{k_3(4)} \\
  &&&& X(2)  \ar[ur]_{X(b)} \ar[ddd]_{k(2)}\\
  &&& A(1) \ar@.[drr]^{\quad \quad A_3(a)} \ar[ddd]_{h(1)}\\
  &~A_3:\ar[ddd]_{h_{3}}&&&& A(3)\oplus D^3_2 \ar@.[r]^{A_3(c)} \ar[ddd]_{h_3(3)} &  A(4)\oplus D^4_3\ar[ddd]_{h_3(4)} \\
  &&&& A(2)  \ar@.[ur]_{A_3(b)}\ar[ddd]_{h(2)}\\
  &&& B(1) \ar@.[drr]^{\quad \quad B_3(a)} \ar[ddd]\\
  &~B_3:\ar[ddd]&&&& B(3)\oplus D^3_2 \ar@.[r]^{B_3(c)} \ar[ddd] &  B(4)\oplus D^4_3\ar[ddd]\\
  &&&& B(2)  \ar@.[ur]_{B_3(b)}\ar[ddd]\\
  &&&0  \\
  &0&&&& 0  &  0 \\
  &&&&0
  }$$
such that $A_{3}\in \Rep{Q}{\widetilde{\sf{F}}}$ and $B_{3}\in\Phi(\widetilde{\sf{C}})$.
This implies that $X$ is in $\mathsf{T}$, and hence completes the proof.\qed
\end{bfhpg}

\bigskip
\section*{Acknowledgments}%%%%%%%%%%%%%%%%%%%%%%%%%%%%%%%%%%%%%%%%%%%%
\noindent
We thank Henrik Holm, Zhongkui Liu, Gang Yang and Xiaoxiang Zhang for helpful discussions related to this work.

\bibliographystyle{amsplain-nodash}

%\bibliographystyle{amsplain}
%\bibliography{../+references}

\def\cprime{$'$}
  \providecommand{\arxiv}[2][AC]{\mbox{\href{http://arxiv.org/abs/#2}{\sf
  arXiv:#2 [math.#1]}}}
  \providecommand{\oldarxiv}[2][AC]{\mbox{\href{http://arxiv.org/abs/math/#2}{\sf
  arXiv:math/#2
  [math.#1]}}}\providecommand{\MR}[1]{\mbox{\href{http://www.ams.org/mathscinet-getitem?mr=#1}{#1}}}
  \renewcommand{\MR}[1]{\mbox{\href{http://www.ams.org/mathscinet-getitem?mr=#1}{#1}}}
\providecommand{\bysame}{\leavevmode\hbox to3em{\hrulefill}\thinspace}
\providecommand{\MR}{\relax\ifhmode\unskip\space\fi MR }
% \MRhref is called by the amsart/book/proc definition of \MR.
\providecommand{\MRhref}[2]{%
  \href{http://www.ams.org/mathscinet-getitem?mr=#1}{#2}
}
\providecommand{\href}[2]{#2}

\end{document}